\documentclass[journal,twoside,web]{IEEEtran}
\usepackage{tmi}
\usepackage{amsmath,amssymb,amsfonts}
\usepackage{algorithmic}
\usepackage{graphicx}
\usepackage{textcomp}
\usepackage{authblk}
\usepackage{mathtools}
\usepackage{csquotes}
\usepackage{moresize}
\usepackage{microtype}
\usepackage[linesnumbered,ruled,noend]{algorithm2e}
\usepackage{tikz}
\usepackage{tikz-3dplot}
\usetikzlibrary{fpu}
\usepackage{xcolor}
\usetikzlibrary{calc,spy,intersections,arrows.meta,angles,quotes,3d,math,perspective}
\usepackage{pgfplots}
\usepackage{epstopdf}
\tikzset{
    myspy/.style={
        spy using outlines={%
            rectangle,
            magnification=4,
            width=1.7cm,
            height=0.7cm,
            connect spies,
            ultra thick,
            draw=red,
            inner sep=0,
            outer sep=0,
        }
    }
}
\definecolor{vertexray}{RGB}{150,40,40}
\definecolor{plotred}{RGB}{220, 50, 47}
\definecolor{plotblue}{RGB}{38, 139, 210}
\definecolor{plotgreen}{RGB}{133, 153, 0}
\def\wda{1}
\def\wdb{2}
\def\wdc{4}

\def\xraytube{%
  \begin{scope}[scale=.1, opacity=1, transform shape]
    \draw (-\wda,0)  --++(-90:\wdc)
       --++(180:\wdb) --++(-90:\wdb) --++(0:\wdb)   --++(0:\wda) --++(0:\wda)
       --++(0:\wdb)   --++(90:\wdb)  --++(180:\wdb) --++(90:\wdc)  --cycle;
\end{scope}
}%

\newcommand{\drawcubeTD}[5][]{%
    \pgfmathsetmacro{\xA}{#2}
    \pgfmathsetmacro{\yA}{#3}
    \pgfmathsetmacro{\zA}{#4}
    \pgfmathsetmacro{\L}{#5}

    \coordinate (#1A)  at (\xA,\yA,\zA);
    \coordinate (#1B)  at (\xA+\L,\yA,\zA);
    \coordinate (#1C)  at (\xA+\L,\yA+\L,\zA);
    \coordinate (#1D)  at (\xA,\yA+\L,\zA);
    \coordinate (#1E)  at (\xA,\yA,\zA+\L);
    \coordinate (#1F)  at (\xA+\L,\yA,\zA+\L);
    \coordinate (#1G)  at (\xA+\L,\yA+\L,\zA+\L);
    \coordinate (#1H)  at (\xA,\yA+\L,\zA+\L);
    \draw[gray!60,#1] (#1D) -- (#1H);
    \draw[gray!60,#1] (#1D) -- (#1C);
    \draw[gray!60,#1] (#1D) -- (#1A);
    \draw[thick, #1] (#1A) -- (#1B) -- (#1C);
    \draw[thick,#1] (#1A) -- (#1E) -- (#1F) -- (#1B);
    \draw[thick,#1] (#1F) -- (#1G) -- (#1C);
    \draw[thick,#1] (#1G) -- (#1H) -- (#1E);
}

\newcommand{\opaqueconnect}[3][]{
        \draw [blue] (#2) -- ($(#3)!2!(#2)$);
        \draw [blue] (#3) -- ($(#2)!2!(#3)$);
        \draw [blue, opacity=0.5] (#2) -- (#3);
}

\usepackage{subcaption}
\usepackage[backend=biber,doi=false,url=false]{biblatex}
\AtBeginBibliography{\small}
\usepackage[urlcolor=black,colorlinks,citecolor={blue!50!black}]{hyperref}
\usepackage{cleveref}
\crefname{equation}{}{}
\Crefname{equation}{}{}
\creflabelformat{equation}{#2(#1)#3}

\crefname{figure}{Figure}{Figures}
\Crefname{figure}{Figure}{Figures}

\crefname{table}{Table}{Tables}
\Crefname{table}{Table}{Tables}
\usepackage{siunitx}
\usepackage{glossaries-extra}
\newabbreviation{ula}{ULA}{unadjusted Langevin algorithm}
\newabbreviation{fbp}{FBP}{filtered backprojection}
\newabbreviation{ct}{CT}{computed tomography}
\newabbreviation{art}{ART}{algebraic reconstruction techniques}
\newabbreviation{sm}{SM}{system matrix}
\newabbreviation{cg}{CG}{conjugate gradient}
\newabbreviation{nag}{NAG}{Nesterov's accelerated gradient descent algorithm}

\usepackage{bm}
\newcommand{\R}{\mathbb{R}}

\newcommand{\origin}{O}

\newcommand{\Rc}{\mathcal{R}}

\newcommand{\dd}{\mathrm{d}}
\renewcommand{\phi}{\varphi}

\newcommand{\ie}{\textit{i.e.}}
\newcommand{\eg}{\textit{e.g.}}
\newcommand{\cf}{\textit{cf.}}

\newcommand{\dvec}{\mathbf{d}}

\newcommand{\vol}{\mathrm{Vol}}
\newcommand{\SystemMatrix}{\mathbf{W}}
\newcommand{\Data}{\mathbf{p}}

\renewcommand{\epsilon}{\varepsilon}
\newcommand{\NumDetectors}{N_d}
\DeclareMathOperator*{\argmin}{arg\,min}
\glsdisablehyper%
\title{Computation of a Consistent System Matrix for Cone-beam Computed Tomography}
\author{
    Josef Simbrunner\thanks{J. S. is with the Division of Neuroradiology, Medical University Graz, Graz, Austria (e-mail: \href{mailto:josef.simbrunner@medunigraz.at}{josef.simbrunner@medunigraz.at}).},
    Clemens Krenn\thanks{C. K. and A. H. are with the Institute of Visual Computing, Graz University of Technology, Graz, Austria (e-mails: \href{mailto:clemens.krenn@student.tugraz.at}{clemens.krenn@student.tugraz.at} and \href{mailto:andreas.habring@tugraz.at}{andreas.habring@tugraz.at}).},
    Martin Zach\thanks{M. Z. is with the Biomedical Imaging Group and the Center for Biomedical Imaging, École Polytechnique Fédérale de Lausanne, Lausanne, Switzerland (e-mail: \href{mailto:martin.zach@epfl.ch}{martin.zach@epfl.ch}).},
    \thanks{Funded in whole or in part by the Austrian Science Fund (FWF) (\href{https://www.fwf.ac.at/en/research-radar/10.55776/COE12}{10.55776/COE12}).}
    Andreas Habring
}

\begin{document}
\maketitle
\begin{abstract}
    We propose a method for the computation of a consistent \gls{sm} for two- and three-dimensional cone-beam \gls{ct}.
    The method relies on the decomposition of the cone-voxel intersection volumes into subvolumes that contribute to distinct detector elements and whose contributions to the system matrix admit exact formulae that can be evaluated without the invocation of costly iterative subroutines.
    We demonstrate that the reconstructions obtained when using the proposed system matrix are superior to those obtained when using common line-based integration approaches with numerical experiments on synthetic and real \gls{ct} data. Moreover, we provide a CUDA implementation of the proposed method.
\end{abstract}

\section{Introduction}
\IEEEPARstart{X}{-ray} \glsfirst{ct} is a non-invasive imaging modality with widespread applications in medicine and engineering.
The goal is to determine the \emph{linear attenuation coefficient} $\mu : \R^d \rightarrow \R_{\geq 0}$, where $d \in \{ 2,3 \}$, of some object of interest from many measurements of the difference in X-ray intensity before entering and after exiting the object.
As in many realistic scenarios in medicine and engineering, we assume that the object is supported on some compact $\Omega \subset \R^d$, and that \( \mu \) is zero outside that domain.

Suppose for the moment that \( d = 2 \), and let \( \mathbf{l}_{\theta, r} : \mathbb{R} \to \mathbb{R}^2, t \mapsto r\bm{\theta}^\perp + t\bm{\theta} \) be a line on the plane \( \mathbb{R}^2\), where \( \bm{\theta} = (\cos(\theta), \sin(\theta)) \) is the unit vector that forms an angle \( \theta \in \mathbb{R}\) with the horizontal axis and \( r \in \mathbb{R} \) is the distance from the origin.
According to Beer--Lambert's law, the intensity \( I \) of a monochromatic X-ray traveling along \(\mathbf{l}_{\theta, r}\) decays as $I^\prime(t) = -\mu(\mathbf{l}_{\theta, r}(t)) I(t)$.
Integrating leads to $-\log\left(I(t_1)/I(t_0)\right) = \Rc\{ \mu \}(\theta, r)$, where \( t_0 \) and \( t_1 \) are real numbers such that all intersections of the line with \( \Omega \) happen at values of \( t \) that fulfill \( t_0 < t < t_1 \), and with the \emph{X-ray transform}~\cite{buzug2008}
\begin{equation}\label{eq:xray_trafo}
    \Rc\{ \mu \}(\theta, r) = \int_{\R} \mu(\mathbf{l}_{\theta, r}(t))\,\dd t.
\end{equation}
% TODO: Somehow incorporate something like the following? Would maybe make nice ties to the continuum inversion vs discrete inversion discussion below.
% In the continuum, Radon has established that \( \mu \) can be reconstructed from R(mu). 

In practice, the measurements are a finite collection of intensity measurements obtained by X-ray detectors with some relation to \cref{eq:xray_trafo}.
Often, this relation is formulated as a \emph{sampling} of the X-ray transform.
% Suppose that a detector is located at some (perpendicular) distance \( r_l \in \R \) to a \emph{line} that passes through the origin and forms an angle \( \theta_k \) with the horizontal axis.
% Then, the measurement obtained by that detector, denoted \( p_{\theta_k, r_l} \), is modeled as \( p_{\theta_k,r_l} = \Rc(\mu)(\theta_k, r_l) + \epsilon \), where \( \epsilon \) models measurement noise (for simplicity, modeled as additive although other models are common in low-photon settings).
In that case, the measurements \( (p_{\theta_k, r_l})_{k,l} \) are modeled as \( p_{\theta_k,r_l} = \Rc\{\mu\}(\theta_k, r_l) + \epsilon \), where \( \epsilon \) models measurement noise (here additive for simplicity, although other models are common in low photon settings).
% TODO the following sentence is true, but it neglects to mention that it would be fine if the ray itself was infinitely thin, which is also an unrealistic assumption even in non-cone-beam systems
This model, however, neglects the fact that detectors in practical systems have an active area over which incoming radiation intensity is integrated.
For example, in a parallel-beam system with detectors that have a width \( w >0 \), a more appropriate model of the relation would be that \begin{equation}
    p_{\theta_k,r_l} = \int_{-w/2}^{w/2} \Rc\{\mu\}(\theta_k, r_l + \delta) \, \dd\delta  + \epsilon.
\end{equation}
% passing through the \emph{strip} obtained by widening the line to an appropriate width.
% TODO: This is a simplification for the sake of argument because its not really a strip with a point source and a fan beam. this should probably still be made clear. Alternatively we could define the radon transform with cone beam in mind, but that sounds even more complicated.
% TODO it should not be too hard to do the proper model for fan beam already here. instead of sweeping r, we should sweep theta (and r i think) over the range defined by the detector edges
It is well established in the literature (\cf~\cite[Section 6]{buzug2008}) that such area-based measurement models are superior.

Most X-ray machines used in clinical practice are equipped with a (one-dimensional or two-dimensional) detector array that consists of \( \NumDetectors \) detectors.
During an acquisition, those detectors are simultaneously illuminated and sampled at a fixed rate as the gantry rotates around the object.
Denoting the total number of sampled timesteps, which corresponds to the different angles of illumination, as $N_p$, the total number of measurements $M$ can therefore be expressed as $M=N_dN_p$.
The measurement can usually be expressed as
\begin{equation}
    \mathbf{p} = \mathcal{W}\{\mu\} + \bm{\epsilon}
    \label{eq:inverse problem}
\end{equation}
where \( \mathcal{W} \) is some linear operator whose image space is \( \R^M \), and \( \bm{\epsilon} \in \R^M \) is measurement noise. 
The image reconstruction problem in \Gls{ct} imaging can, therefore, be reduced to the problem of finding some $\mu$ that is consistent with \cref{eq:inverse problem}. 

Generally, we can distinguish two different approaches for solving \cref{eq:inverse problem} depending on the order of inversion and discretization: As an operator between function spaces, $\Rc$ admits an inverse referred to as the \gls{fbp}. The first approach therefore consists of applying a discretization of the \gls{fbp} $\Rc^{-1}$ to the measurement $p$. Due to ill-posedness, however, the \gls{fbp} suffers from instability, especially for few measurements~\cite[Chapter 4]{natterer1986}. The second approach, is to first discretize $\Rc$ in the \emph{forward direction} and afterwards solve the corresponding finite dimensional problem.
The latter approach is the one we focus on in this work as it offers several key benefits, in particular, (i) it is amenable to regularization techniques which are crucial to counter ill-posedness, (ii) it easily allows for adoption to specific geometries, (iii) it allows to properly take into account finite detector widths and different detector sensitivities, and (iv) rays/beams running through objects that potentially produce inconsistencies in the Radon space can be weighted appropriately~\cite{buzug2008}.
%TODO incorporate properly
In particular, in this work we propose a new method for the computation of a consistent discretization of $\mathcal{W}$, that is, the computation of a \glsfirst{sm} $\mathbf{W} \in\R^{M\times N}$ such that we can approximate \eqref{eq:inverse problem} as
\begin{equation}\label{eq:ip}
    \mathbf{p} = \mathbf{W} \mathbf{u}+\bm{\epsilon},
\end{equation}
where $\mathbf{u}\in \R^N$ is a discretization of $\mu$ as a piecewise constant function on a pixel- or voxel grid.

\begin{figure}
        \captionsetup[subfigure]{font=footnotesize}
        \subcaptionbox{Line-based discretization.\label{fig:ct_discretization_line}}[.25\textwidth]{%
        \resizebox{.23\textwidth}{!}{%
        \begin{tikzpicture}
            % Draw a grid of pixels
            \foreach \x in {0,1,2,3} {
                \foreach \y in {0,1,2,3} {
                    \draw[gray] (\x,\y) rectangle +(1,1);
                }
            }
            
            % Draw a straight line crossing through the grid
            \newcommand{\slope}{.5}
            \draw[red, thick] (-0.5,1) -- (-.5+6,1+6*\slope);
            \newcommand{\xc}{-.5+6}
            \newcommand{\yc}{1+6*\slope}
            \draw[fill=green!30] (\xc-.1,\yc-.5) rectangle (\xc+.1,\yc+.5);
            
            % Mark two consecutive cross points
            \newcommand{\x}{2.5}
            \newcommand{\xs}{-0.5 + \x}
            \newcommand{\ys}{1+\x*\slope}
            \fill[blue] (\xs,\ys) circle (0.1);

            \newcommand{\xx}{3.5}
            \newcommand{\xe}{-0.5 + \xx}
            \newcommand{\ye}{1+\xx*\slope}
            \fill[blue] (\xe,\ye) circle (0.1);

            \draw[blue, thick,<->] (\xs,\ys+.2) -- (\xe,\ye+.2);

            % Label the line segment between cross points
            \newcommand{\xl}{0.5*\xs+0.5*\xe}
            \newcommand{\yl}{0.3+0.5*\ys+0.5*\ye}
            \node[blue] at (2.5,3.2) {$w_{i,j}$};

            % emphasize the relevant pixel
            \draw[thick] (2,2) rectangle +(1,1);
            \node at (2.6,1.75) {voxel $j$};

            \node[red] at (0.4,1.75) {X-ray};
            \node at (4.6,4.2) {detector $i$};
        \end{tikzpicture}}}%
        \subcaptionbox{Area-based discretization.\label{fig:ct_discretization_strip}}[.25\textwidth]{
            \resizebox{.23\textwidth}{!}{%
            \begin{tikzpicture}
                % Draw a grid of pixels
                \foreach \x in {0,1,2,3} {
                    \foreach \y in {0,1,2,3} {
                        \draw[gray] (\x,\y) rectangle +(1,1);
                    }
                }
                
                % Fill the area
                \newcommand{\slope}{.6}
                \newcommand{\slopel}{.4}
                \fill[blue!30,opacity=.8] (-.5+2.5,1+2.5*\slopel) -- (-.5+3.5,1+3.5*\slopel) -- (3,3) -- (-.5+3.333,3) --(-.5+2.5,1+2.5*\slope) -- cycle;

                % Draw a straight line crossing through the grid
                \draw[red, thick] (-0.5,1) -- (-.5+6,1+6*\slope);
                \draw[red, thick] (-0.5,1) -- (-.5+6,1+6*\slopel);
                \newcommand{\xl}{-.5+6}
                \newcommand{\yu}{1+6*\slope}
                \newcommand{\yl}{1+6*\slopel}
                \draw[fill=green!30] (\xl,\yl) rectangle (\xl+.2,\yu);

                % Label the line segment between cross points
                \renewcommand{\xl}{0.5*\xs+0.5*\xe}
                \renewcommand{\yl}{0.3+0.5*\ys+0.5*\ye}
                \node[blue] at (2.5,2.4) {$w_{i,j}$};
    
                % emphasize the relevant pixel
                \draw[thick] (2,2) rectangle +(1,1);
                % \node at (2.6,1.75) {pixel $j$};
    
            \end{tikzpicture}}
        }
    \caption{Discretization of the \gls{ct} measurement process for the fan beam geometry.}
    \label{fig:ct_discretization}
\end{figure}
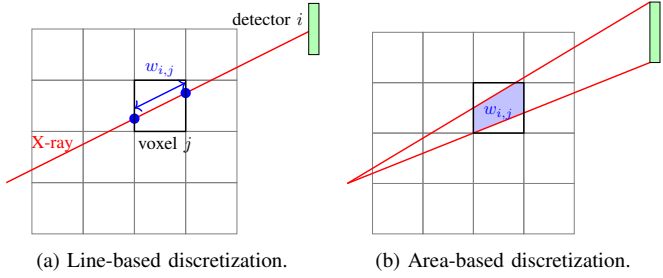

The term \emph{consistent} refers to the fact that the discretized projection is exact for attenuation coefficients $\mu$ which are, indeed, piecewise constant on the pixel/voxel grid, \ie, in this case we have
\[
    \mathcal{W}\{\mu\} = \SystemMatrix \mathbf{u}.
\]
Due to the computational difficulty, in practice $\SystemMatrix$ is usually approximated by computing $w_{i,j}$ as the path length within voxel $j$ of a single ray to detector element $i$.
Such a line-based integration, however, is inconsistent in the above introduced sense.
In \cref{fig:ct_discretization} we illustrate line-based and consistent, volume-based integration in the 2D case.

\subsection{Contributions} We provide the following contributions;
\begin{enumerate}
    \item We derive explicit formulas for the computation of a consistent, area-/volume-based \gls{sm} for cone-beam flat-detector \gls{ct} in two and three dimensions.
    \item We experimentally verify correctness of our formulas as well as practical advantages of the consistent \gls{sm} in certain resolution and/or regularization regimes.
    \item We provide a CUDA implementation of the proposed method available at~\cite{volumebased_code}.
\end{enumerate}
\section{Related work}
Most existing methods for the computation of the system matrix, such as those described in~\cite{jacobs1998fast,joseph1982improved,siddon1985fast,wang2021fast,zhang2014fast,zhang2023iterative,zhao2003fast}, rely on line-based integration due to the computational simplicity.
In~\cite{lo1988strip}, Lo demonstrated the advantages of using exact weighting coefficients in medical image reconstruction for parallel beam projections, where exact weighting factors can be easily calculated. 
In \cite{sen1995exact,ha2015study} for 3D cone beam \gls{ct} the \gls{sm} is computed by successively clipping the voxels by intersecting them with faces of the X-ray bundle and afterward computing the sub-voxel volume using triangulation approaches.
In \cite{sen1995exact} the authors found that their exact weighting factors are superior to approximate ones in medical image reconstructions. 
In \cite{yao2009analytically}, the authors express the weighting coefficients as a three-dimensional integral, which they reduce to a one-dimensional one that distinguishes various cases.
They numerically demonstrate the advantage of the volume-based weighting factors over line-based ones.
The authors of \cite{ha2018look} circumvent heavy on-line computations through the pre-computation of look-up tables generated by polynomial approximation of the weighting coefficients.

To the best of our knowledge all available methods rely on costly subroutines (triangulation, clipping, numerical integration). To the contrary, we derive explicit formulae for computing the weighting factors. Moreover, contrary to the cited works we publicly provide a CUDA implementation of the method.

\section{Method}
\label{method}
In the following we will present our approach to explicitly and analytically compute $\SystemMatrix$.
Some of the more tedious computations are postponed to the appendix~\cref{sec:appendix}.

\subsection{Notation and Setup}

We will rely on two fundamental geometrical objects: \emph{points} and \emph{rays}.
We denote \emph{points} in $\R^d$ with italic capital letters and denote the Euclidean distance between two points $P,Q\in\R^d$ as $|PQ|$.
We denote \emph{rays} in \( \R^d \) with bold roman lowercase letters.
A ray is understood as the unidirectional infinite extension of a line segment that connects two points.

We now describe our assumptions on the geometrical setup.
Here, we stick to a three-dimensional description for simplicity; the two-dimensional case is obtained by dropping the third coordinate.
We denote the origin of the coordinate system as $\origin$. An X-ray source (idealized as a point-source) is initially located at \( S = (s, 0, 0) \), where \( s > 0 \) is the distance from the origin.
The system is equipped with a flat-panel detector array with rectangular detectors that are characterized by their sidelenghts \( d_y, d_z > 0 \) in the \( y \) and \( z \) direction, respectively, and whose pixel response is uniform.
This flat-panel detector array is oriented perpendicular to the ray that connects the X-ray point-source with the origin and located at a distance \( d \) from the origin.
We index the corners of the detector grid with a vector \( \mathbf{m} = (\mathbf{m}_y, \mathbf{m}_z) \in \mathbb{Z}^2 \) so that the edge identified by \( \mathbf{m} \) is initially located at the point \( (-d, \mathbf{m}_yd_y, \mathbf{m}_zd_z) \).
To make our setup more concrete, we specify that after a rotation of the source and the detector array by an angle \( \varphi \) around the \( z \) axis, the ray from \( S \) through the detector edge identified by \( \mathbf{m} \) is given by
\begin{equation}
	\dvec_{\mathbf{m}} = \begin{pmatrix}
			s\cos\varphi \\
			s\sin\varphi \\
			0
			\end{pmatrix} + \lambda\begin{pmatrix}
			(s + d)\cos\varphi + \mathbf{m}_{y}d_{y}\sin\varphi \\
			(s + d)\sin\varphi - \mathbf{m}_{y}d_{y}\cos\varphi \\
			- \mathbf{m}_{z}d_{z}
	\end{pmatrix}.
 \label{eq:z rotation}
\end{equation}
We denote the angle that is enclosed by the rays from the source through the origin and from the source through the detector edge identified by \( \mathbf{m} \) from the top-view as \( \beta_{\mathbf{m}_y} \) (see, \eg, \(\beta_{-3}\) in~\cref{fig:2d triangle}).

We consider orthorhombic voxels with the sidelengths \( a \), \( b \) and \( c \) in \( x \), \( y \), and \( z \) direction, respectively and the origin of the coordinate system is centered in one of those voxels.
We index this voxel grid with \( \mathbf{n} = (\mathbf{n}_x, \mathbf{n}_y, \mathbf{n}_z) \in \mathbb{Z}^3 \) so that the voxel identified by \( \mathbf{n} \) has the vertices
\begin{align}
    V_{\mathbf{n},01} &= ({\mathbf{n}_{x}}a, {\mathbf{n}_{y}b}, {\mathbf{n}_{z}c}), \\
    V_{\mathbf{n},02} &= (({\mathbf{n}_{x}+1)a}, {\mathbf{n}_{y}b}, {\mathbf{n}_{z}c}), \\
    V_{\mathbf{n},03} &= (({\mathbf{n}_{x}+1)a}, ({\mathbf{n}_{y}} + 1){b}, {\mathbf{n}_{z}c}),\text{ and}\\
    V_{\mathbf{n},04} &= ({\mathbf{n}_{x}}{a}, ({\mathbf{n}_{y}}+1){b}, {\mathbf{n}_{z}c})
\end{align}
on its bottom face and the vertices 
\begin{align}
    V_{\mathbf{n},11} &= ({\mathbf{n}_{x}}{a}, {\mathbf{n}_{y}b}, ({\mathbf{n}_{z}} + 1){c}), \\
    V_{\mathbf{n},12} &= (({\mathbf{n}_{x}+1)a}, {\mathbf{n}_{y}b}, ({\mathbf{n}_{z}} + 1){c}), \\
    V_{\mathbf{n},13} &= (({\mathbf{n}_{x}+1)a}, ({\mathbf{n}_{y}} + 1){b}, ({\mathbf{n}_{z}} + 1){c}),\text{ and}\\
    V_{\mathbf{n},14} &= ({\mathbf{n}_{x}}{a}, ({\mathbf{n}_{y}} + 1){b}, ({\mathbf{n}_{z}} + 1){c}).
\end{align}
We denote the angle from the top view that is enclosed by the rays from the source through the origin and from the source through the vertex $V_{\mathbf{n},ij}$ (which depends on the gantry rotation angle \( \varphi \)) 
as \( \gamma_{\mathbf{n},j,\phi} \) (see again, \eg, \( \gamma_{\mathbf{n},2,\phi}\) in~\cref{fig:2d triangle}).

\subsection{The Two-dimensional Case}
We begin with the two-dimensional case where $z\equiv 0$.
The results that we obtain in this section will also be useful for the three-dimensional case.
We denote the vertices of the voxel identified by the two-dimensional index \( \mathbf{n} = (\mathbf{n}_x, \mathbf{n}_y) \) as \( V_{\mathbf{n},1}, V_{\mathbf{n},2}, V_{\mathbf{n},3} \) and \( V_{\mathbf{n},4} \).
Though one integer suffices to index the detectors in this case, we still denote this integer with as an upright, bold \( \mathbf{m} = \mathbf{m}_y \).

We can assume without loss of generality that the angles to the vertices verify \( \max(\gamma_{\mathbf{n},1,\phi}, \gamma_{\mathbf{n},2,\phi}) \leq \min(\gamma_{\mathbf{n},3,\phi}, \gamma_{\mathbf{n},4,\phi}) \).
In the converse case we may rotate the coordinate system by \(\pi/2\) until this condition is satisfied.
The index \( \mathbf{m} \) of any detector element that intersects the cone spanned by the source and the voxel satisfies
\begin{equation}
\begin{split}
	\tan\left( \min(\gamma_{\mathbf{n},1,\phi}, \gamma_{\mathbf{n},2,\phi}) \right)\frac{s + d}{d_{y}} < \mathbf{m} &\\
    < \tan\left( \max(\gamma_{\mathbf{n},3,\phi}, \gamma_{\mathbf{n},4,\phi}) \right)&\frac{s + d}{d_{y}}.
    \label{eq:relevant detectors 2d}
\end{split}
\end{equation}
The proposed algorithm calculates the contributions of the voxel to each detector element by distinguishing between different cases of the relationships of the angles \( \gamma_{\mathbf{n},i,\phi}\), \(i=1,2,3,4\) and \( \beta_{\mathbf{m}} \) for all \( \mathbf{m} \) satisfying \cref{eq:relevant detectors 2d}, \cf~\cref{fig:2d triangle}.

\begin{figure}
    \centering
    \begin{tikzpicture}[line cap=round,line join=round,>=Latex,scale=1]
        \pgfmathsetmacro{\anglePhi}{-69};
        \pgfmathsetmacro{\dS}{5};
        \pgfmathsetmacro{\Sx}{\dS * cos(\anglePhi)};
        \pgfmathsetmacro{\Sy}{\dS * sin(\anglePhi)};
        \pgfmathsetmacro{\dD}{5};
        \pgfmathsetmacro{\Dx}{-\dD * cos(\anglePhi)};
        \pgfmathsetmacro{\Dy}{-\dD * sin(\anglePhi)};
        \pgfmathsetmacro{\ddy}{0.8};
        \pgfmathsetmacro{\sidelength}{1.8};
        \pgfmathsetmacro{\Vxx}{\sidelength / 2}
        
        \coordinate (Vone) at (-\Vxx-\sidelength, -\Vxx+\sidelength);
        \coordinate (Vtwo) at (-\Vxx, -\Vxx+\sidelength);
        \coordinate (Vfour) at (-\Vxx-\sidelength, \Vxx+\sidelength);
        \draw[step=\sidelength, help lines,shift={(\Vxx,\Vxx)}] (-4cm,-4cm) grid (2.2cm, 2.2cm);
        \node [below left] at (Vone) {\( V_{\mathbf{n},2} \)};
        \node [above left] at (Vfour)  {\( V_{\mathbf{n},1} \)};

        \begin{scope}[shift={(\Sx,\Sy)},rotate=120-67-15]
            \xraytube
        \end{scope}
        \draw (-\Vxx-\sidelength,-\Vxx+\sidelength) rectangle ++(\sidelength, \sidelength);
        \coordinate (O) at (0, 0);
        \coordinate (E1) at (0, -1);
        \coordinate (S) at (\Sx, \Sy);
        \coordinate (D0) at (\Dx, \Dy);
        \draw [->, opacity=0.2] (S) -- (D0);
        \foreach \di in {0,1,2,3}
        {
            \pgfmathtruncatemacro{\diPlusOne}{\di + 1};
            \draw (D\di) -- ++({-\ddy*sin(\anglePhi)}, {+\ddy*cos(\anglePhi)}) coordinate (D\diPlusOne);
            \draw (D\diPlusOne) -- ++({-\ddy/2*cos(\anglePhi)}, {-\ddy/2*sin(\anglePhi)}) node [left,rotate=\anglePhi] {\( \diPlusOne\)};
            % \draw [->, opacity=0.2] (S) -- (D\diPlusOne);
        }
        \draw (D0) -- ++({-\ddy/2*cos(\anglePhi)}, {-\ddy/2*sin(\anglePhi)}) node [left,rotate=\anglePhi] {\( i = 0 \)};
        \foreach \di in {0,-1,-2,-3}
        {
            \pgfmathtruncatemacro{\diPlusOne}{\di - 1};
            \draw (D\di) -- ++({+\ddy*sin(\anglePhi)}, {-\ddy*cos(\anglePhi)}) coordinate (D\diPlusOne);
            \draw (D\diPlusOne) -- ++({-\ddy/2*cos(\anglePhi)}, {-\ddy/2*sin(\anglePhi)}) node [left,rotate=\anglePhi] {\( \diPlusOne\)};
            \draw [->, opacity=0.2] (S) -- (D\diPlusOne);
        }
        \draw [->] (S) -- (D-3);
        \draw [dashed] (D4) -- ++ ({-\ddy*sin(\anglePhi)}, {\ddy*cos(\anglePhi)});
        \draw [dashed] (D-4) -- ++ ({\ddy*sin(\anglePhi)}, {-\ddy*cos(\anglePhi)});
        \draw [->] (0,0) -- ++(0,-1) node[right]{\( x \)};
        \draw [->] (0,0) -- ++(1,0) node[above]{\( y \)};
        \foreach \vertex in {Vone, Vfour} {
            \draw [->, vertexray] (S) -- (\vertex);
        }
        \pic [draw, ->, "$\phi$", angle eccentricity=1.5] {angle = E1--O--S};
        \path [name path=v12] (Vone) -- (Vtwo);
        \path [name path=sd3] (S) -- (D-3);
        \path [overlay, name intersections={of=v12 and sd3, by=px}];
        \fill (px) circle (0.05cm) node [above right] {\( P_{\mathbf{m},\mathbf{n},\phi,x} \)};

        \path [name path=v14] (Vone) -- (Vfour);
        \path [overlay, name intersections={of=v14 and sd3, by=py}];
        \fill (py) circle (0.05cm) node [left] {\( P_{\mathbf{m},\mathbf{n},\phi,y} \)};
        \pic [draw, ->, "$\gamma_{\mathbf{n},2,\phi}$", angle eccentricity=1.1, angle radius=140] {angle = O--S--Vone};
        \pic [draw, ->, "$\beta_{-3}$", angle eccentricity=1.1, angle radius=100] {angle = O--S--D-3};
        \fill [blue!50!black, opacity=0.5] (px) -- (py) -- (Vone) -- cycle;

        % \node at ($(S)!0.5!(Vone)$) [rotate=-40, below, vertexray] {\( \mathbf{g}_1 \)};
        \node at ($(S)!0.9!(D-3)$) [rotate=-50, below] {\( \mathbf{b}_{-3} \)};
    \end{tikzpicture}
    \caption{%
        Schematic subdivision of a pixel by rays to the detector edges and/or pixel vertices. We denote $\beta_{\mathbf{m}}$ the angles to the rays to detector edges and $\gamma_{\mathbf{n},i,\phi}$ the angles of rays to vertices of the pixel. 
        % \( g_1 \), \( g_2 \), \( g_3 \), and \( g_4 \) (green) through the vertices \( V_1 \), \( V_2 \), \( V_3 \), and \( V_4 \), respectively.
        % The rays \( g_1 \) and \( b_i \) delimit the triangle \( V_1 P_X P_Y \).
        % Also shown are their respective angles \( \beta_g \) and \( \beta_i \).
    }%
    \label{fig:2d triangle}
\end{figure}
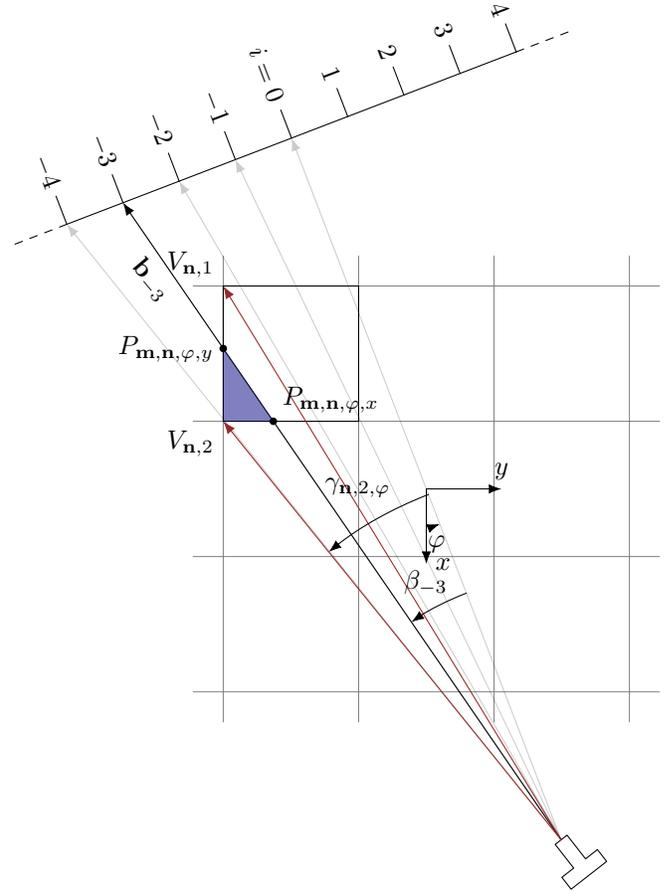

\begin{figure}
	\centering
    \begin{tikzpicture}[line cap=round,line join=round,>=Latex,scale=1]
        \pgfmathsetmacro{\anglePhi}{-83};
        \pgfmathsetmacro{\dS}{5};
        \pgfmathsetmacro{\Sx}{\dS * cos(\anglePhi)};
        \pgfmathsetmacro{\Sy}{\dS * sin(\anglePhi)};
        \pgfmathsetmacro{\dD}{5};
        \pgfmathsetmacro{\Dx}{-\dD * cos(\anglePhi)};
        \pgfmathsetmacro{\Dy}{-\dD * sin(\anglePhi)};
        \pgfmathsetmacro{\ddy}{0.8};
        \pgfmathsetmacro{\sidelength}{1.8};
        \pgfmathsetmacro{\Vxx}{\sidelength / 2}
        
        \coordinate (Vone) at (-\Vxx-\sidelength, -\Vxx+\sidelength);
        \coordinate (Vtwo) at (-\Vxx, -\Vxx+\sidelength);
        \coordinate (Vthree) at (-\Vxx, \Vxx+\sidelength);
        \coordinate (Vfour) at (-\Vxx-\sidelength, \Vxx+\sidelength);
        \draw[step=\sidelength, help lines,shift={(\Vxx,\Vxx)}] (-4cm,-4cm) grid (2.2cm, 2.2cm);
        % \node [right] at (Vtwo) {\( V_{\mathbf{n},2} \)};
        % \node [above left] at (Vfour)  {\( V_{\mathbf{n},1} \)};

        \begin{scope}[shift={(\Sx,\Sy)},rotate=120-67-15]
            \xraytube
        \end{scope}
        \draw (-\Vxx-\sidelength,-\Vxx+\sidelength) rectangle ++(\sidelength, \sidelength);
        \coordinate (O) at (0, 0);
        \coordinate (E1) at (0, -1);
        \coordinate (S) at (\Sx, \Sy);
        \coordinate (D0) at (\Dx, \Dy);
        \draw [->, opacity=0.2] (S) -- (D0);
        \foreach \di in {0,1,2,3}
        {
            \pgfmathtruncatemacro{\diPlusOne}{\di + 1};
            \draw (D\di) -- ++({-\ddy*sin(\anglePhi)}, {+\ddy*cos(\anglePhi)}) coordinate (D\diPlusOne);
            \draw (D\diPlusOne) -- ++({-\ddy/2*cos(\anglePhi)}, {-\ddy/2*sin(\anglePhi)}) node [left,rotate=\anglePhi] {\( \diPlusOne\)};
            % \draw [->, opacity=0.2] (S) -- (D\diPlusOne);
        }
        \draw (D0) -- ++({-\ddy/2*cos(\anglePhi)}, {-\ddy/2*sin(\anglePhi)}) node [left,rotate=\anglePhi] {\( i = 0 \)};
        \foreach \di in {0,-1,-2,-3}
        {
            \pgfmathtruncatemacro{\diPlusOne}{\di - 1};
            \draw (D\di) -- ++({+\ddy*sin(\anglePhi)}, {-\ddy*cos(\anglePhi)}) coordinate (D\diPlusOne);
            \draw (D\diPlusOne) -- ++({-\ddy/2*cos(\anglePhi)}, {-\ddy/2*sin(\anglePhi)}) node [left,rotate=\anglePhi] {\( \diPlusOne\)};
            \draw [->, opacity=0.2] (S) -- (D\diPlusOne);
        }
        \draw [->,name path=sd3] (S) -- (D-3);
        \draw [->,name path=sd2] (S) -- (D-2);
        \draw [dashed] (D4) -- ++ ({-\ddy*sin(\anglePhi)}, {\ddy*cos(\anglePhi)});
        \draw [dashed] (D-4) -- ++ ({\ddy*sin(\anglePhi)}, {-\ddy*cos(\anglePhi)});
        % \draw [->] (0,0) -- ++(0,-1) node[right]{\( x \)};
        % \draw [->] (0,0) -- ++(1,0) node[above]{\( y \)};
        \foreach \vertex in {Vtwo, Vfour} {
            \draw [->, vertexray] (S) -- (\vertex);
        }
        % \pic [draw, ->, "$\phi$", angle eccentricity=1.5] {angle = E1--O--S};
        \path [name path=v12] (Vone) -- (Vtwo);
        \path [name path=v43] (Vfour) -- (Vthree);
        \path [overlay, name intersections={of=v12 and sd3, by=Q1}];
        \path [overlay, name intersections={of=v43 and sd3, by=P1}];
        \fill (Q1) circle (0.05cm) node [below left] {\( P_{\mathbf{m},\mathbf{\tilde{n}},\phi,x} \)};
        \fill (P1) circle (0.05cm) node [above left] {\( P_{\mathbf{m},\mathbf{n},\phi,x} \)};

        \path [overlay, name intersections={of=v12 and sd2, by=Q2}];
        \path [overlay, name intersections={of=v43 and sd2, by=P2}];
        \fill (Q2) circle (0.05cm) node [above right] {\( P_{\mathbf{m}+1,\mathbf{\tilde{n}},\phi,x} \)};
        \fill (P2) circle (0.05cm) node [above right] {\( P_{\mathbf{m}+1,\mathbf{n},\phi,x} \)};

        \pic [draw, ->, "$\beta_{\mathbf{b}_{-3}}$", angle eccentricity=1.1, angle radius=110] {angle = O--S--D-3};
        \pic [draw, ->, "$\beta_{\mathbf{b}_{-2}}$", angle eccentricity=1.1, angle radius=80] {angle = O--S--D-2};
        \fill [blue!50!black, opacity=0.5] (P1) -- (P2) -- (Q2) -- (Q1) -- cycle;

        \node at ($(S)!0.9!(D-3)$) [rotate=-68, below] {\( \mathbf{b}_{-3} \)};
        \node at ($(S)!0.9!(D-2)$) [rotate=-74, below] {\( \mathbf{b}_{-2} \)};
    \end{tikzpicture}
    \caption{
        Trapezoid shape sub-voxel delimited by two rays to subsequent detector elements. 
    }%
    \label{fig:2d trapezoid}
\end{figure}
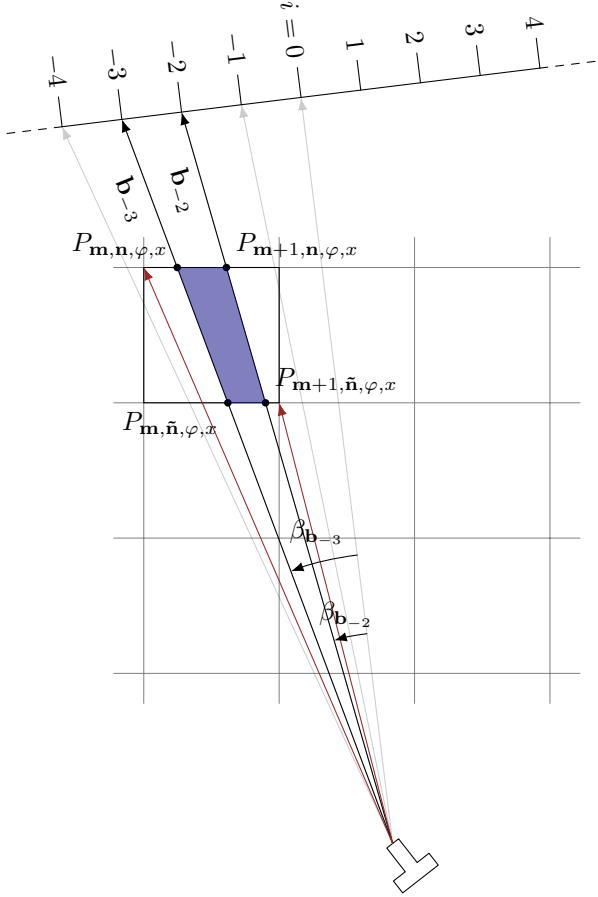
For this purpose, the formulas of two fundamental area elements are sufficient:
a triangle and a trapezoid.
With these fundamental atoms, the area of the intersection of the X-ray cone from the source to any detector element can be calculated either directly or by an appropriate summation or subtraction of such area elements.
\subsubsection{Triangle}
Assume again without loss of generality that one of the corners of the triangle is \( V_{\mathbf{n},1} \) or $V_{\mathbf{n},2}$~\cref{fig:2d triangle}.
There, we specifically depict the case that the ray from the source through the detector identified by \( \mathbf{m} \) intersects the voxel at the points \( P_{\mathbf{m},\mathbf{n},\phi,x}\), and \(P_{\mathbf{m},\mathbf{n},\phi,y} \) such that \( P_{\mathbf{m},\mathbf{n},\phi,x}\), \(P_{\mathbf{m},\mathbf{n},\phi,y} \) and \( V_{\mathbf{n},1} \) form a triangle.
However, the derivation is analogous if $V_{\mathbf{n},2}$ is enclosed by the ray.
Using \eqref{eq:x expression} and \eqref{eq:y expression}, the intersection points are given by
\begin{equation}
\begin{aligned}
    &P_{\mathbf{m},\mathbf{n},\varphi,x} \\
    &\hspace{0.2cm}= \begin{pmatrix}
		(\mathbf{n}_{x} + 1)a \\
		s\sin\varphi - (s\cos\varphi - \left( \mathbf{n}_{x} + 1 \right)a)\tan{(\varphi} - \beta_{\mathbf{m}})
	\end{pmatrix},
\end{aligned}
\end{equation}
and
\begin{equation}
	P_{\mathbf{m},\mathbf{n},\varphi,y} = \begin{pmatrix}
		s\cos\varphi - (s\sin\varphi - \mathbf{n}_{y}b)\cot{(\varphi} - \beta_{\mathbf{m}}) \\
		\mathbf{n}_{y}b
	\end{pmatrix},
\end{equation}
and the resulting area, relative to the area \( ab \) of the voxel, can be expressed solely by the rotation angle of the gantry \( \varphi \), the angles enclosed by the source-origin and source-vertex rays \( \gamma_{\mathbf{n},1,\phi}, \gamma_{\mathbf{n},2,\phi} \), and the angle enclosed by the source-origin and source-detector ray \( \beta_{\mathbf{m}} \) as
\begin{equation}
\begin{split}
	&A_{3}(\varphi, \gamma_{\mathbf{n},2,\phi}, \gamma_{\mathbf{n},1,\phi}, \beta_{\mathbf{m}}) \\
    &\quad= \frac{a}{b\left| \sin{(2\varphi - {2\beta}_{i})} \right|}\left( \frac{\sin{(\varphi} - \gamma_{\mathbf{n},1,\phi})\sin{(\beta_{i}} - \gamma_{\mathbf{n},2,\phi})}{\sin{(\gamma_{\mathbf{n},1,\phi}} - \gamma_{\mathbf{n},2,\phi})} \right)^{2}.
\end{split}
\end{equation}
For the edge cases in which \( \gamma_{\mathbf{n},1,\phi} = \gamma_{\mathbf{n},2,\phi} \) or \( \beta_{\mathbf{m}} = \varphi \), the area factor is zero.
\subsubsection{Trapezoid}%
\label{trapezoid}
% The trapezoid is confined by two rays which intersect two opposite pixel sides.
If the ray from the source through the detector corner identified by \( \mathbf{m}\) and the ray from the source through the detector corner identified by \( \mathbf{m} + 1 \) both enter and exit the voxel identified by \( \mathbf{n} \) on opposite sides, the intersection volume is a trapezoid.
This is illustrated in \cref{fig:2d trapezoid}.
Using the previous notation, the points of intersection can be enumerated as \( P_{\mathbf{m},\mathbf{n},\varphi,x}, P_{\mathbf{m} + 1,\mathbf{n},\varphi,x}, P_{\mathbf{m},\tilde{\mathbf{n}},\varphi,x}\), and \( P_{\mathbf{m} + 1,\tilde{\mathbf{n}},\varphi,x} \),.
where \( \tilde{\mathbf{n}} = \mathbf{n} + (0, 1) \).
As the height of the trapezoid is \( a \), the area of that trapezoid is
\begin{equation}
\begin{split}
    &\tfrac{a}{2}\bigl(\left| P_{\mathbf{m},\mathbf{n},\varphi,x} P_{\mathbf{m} + 1,\mathbf{n},\varphi,x}\right| + \left| P_{\mathbf{m},\tilde{\mathbf{n}},\varphi,x}P_{\mathbf{m} + 1,\tilde{\mathbf{n}},\varphi,x}  \right|\bigr),
 \label{eq:area trapezoid}
 \end{split}
\end{equation}
and the area relative to the area \( ab \) of the voxel can be expressed as
\begin{equation}\label{eq:area trapezoid rel}
\begin{split}
	A_{4}( \varphi, \mathbf{n}_{x},&\beta_{\mathbf{m}},\beta_{\mathbf{m} + 1} ) \\
    =&\tfrac{1}{b}\bigl(\left| (s\cos\varphi - \left( \mathbf{n}_{x} + \tfrac{1}{2} \right)a)\right|\\
    &\times \left| \tan{(\varphi} - \beta_{\mathbf{m}}) - \tan{(\varphi} - \beta_{\mathbf{m} + 1}) \right|\bigr).
\end{split}
\end{equation}
% Using~\cref{eq:cos phi nx a} and~\cref{eq:cos phi nx1 a}, \cref{eq:area trapezoid rel} can be written as:
% \begin{equation}
% \begin{split}
% 	&A_{4\mathrm{rel}}\left( \varphi,\betago,\betagtwo,\beta_{i},\beta_{i + 1} \right) \\ 
%     = &\frac{a}{b}\left| \frac{\sin{(2\varphi} - \betago - \betagtwo)}{2\sin{(\betagtwo} - \betago)}\right|\\
%     &\times\left| \tan{(\varphi} - \beta_{i}) - \tan{(\varphi} - \beta_{i + 1}) \right|
% \end{split}
% \end{equation}
\subsubsection{Area Factors}%
\label{area-factors}
% We now detail the algorithm that we use to combine the atomic area elements to obtain the area of the intersection of the source-detector cone with the voxel.
Building on these formulas, we can compute the entire area of the sub-voxel contributing to detector element $\mathbf{m}$ as follows:
Starting from the left-most corner vertex of the voxel (\ie, the one with minimal angle $\gamma_{\mathbf{n},i,\phi}$), we sweep the rays to the detector elements and the remaining corners of the voxel.
All the emerging sub-voxel areas can be computed as sums or differences of triangles and/or trapezoids
In mathematical notation, let us first denote the sequence $(\gamma_{\mathbf{n}, i, \phi})_{i=1}^4$ sorted in ascending order as $(\Gamma_{\mathbf{n}, i, \phi})_{i=1}^4$.
The lines spanned by the angles \( \Gamma_{\mathbf{n}, 2, \phi} \) and \( \Gamma_{\mathbf{n}, 3, \phi} \) split the pixel into three parts: two triangles which enclose a trapezoid.
The rays \( \beta_{\mathbf{m}} \) can be assigned to these three sub-cells according to their magnitudes so that for the first triangle \( \Gamma_{\mathbf{n}, 1 \phi} \leq \beta_{\mathbf{m}} \leq \Gamma_{\mathbf{n}, 2 \phi} \), for the trapezoid \( \Gamma_{\mathbf{n}, 2 \phi} < \beta_{\mathbf{m}} \leq \Gamma_{\mathbf{n}, 3 \phi} \) and for the second triangle \( \Gamma_{\mathbf{n}, 3 \phi} < \beta_{\mathbf{m}} \leq \Gamma_{\mathbf{n}, 4 \phi} \).

As an example\footnote{ignoring $\phi$ for a moment}, let \( w_{\mathbf{m},\mathbf{n}} \) for \( \mathbf{m} = -4, \dotsc, 0 \) be the area factors for the case depicted in~\cref{fig:2d triangle}.
Then \(f_i\) are given by
\begin{equation}\label{eq:2d procedure}
\begin{aligned}
    \begin{pmatrix}
        w_{-4,\mathbf{n}} \\
        w_{-3,\mathbf{n}}^1 \\
        w_{-3,\mathbf{n}}^2\\
        w_{-2,\mathbf{n}}^1\\
        w_{-2,\mathbf{n}}^2\\
        w_{-1,\mathbf{n}}\\
        w_{0,\mathbf{n}}
    \end{pmatrix} = &\begin{pmatrix}
        1 & 0 & 0 &0 &0 &0 &0 \\
        -1 & 1 & 0 &0 &0 &0 &0 \\
        0 & 0 & 1 &0 &0 &0 &0 \\
        0 & 0 & 0 &1 &0 &0 &0 \\
        0 & 0 & 0 &0 &1 &-1 &0 \\
        0 & 0 & 0 &0 &0 &1 &-1 \\
        0 & 0 & 0 &0 &0 &0 &1 \\
    \end{pmatrix}\\
    &\times\begin{pmatrix}
    		A_{3}(\varphi, \Gamma_{\mathbf{n},1,\phi}, \Gamma_{\mathbf{n},2,\phi}, \beta_{-3})\\
    		A_{3}(\varphi, \Gamma_{\mathbf{n},1,\phi}, \Gamma_{\mathbf{n},2,\phi}, \Gamma_{\mathbf{n},2,\phi})\\
    		A_{4}(\varphi,\mathbf{n}_x, \Gamma_{\mathbf{n},2,\phi}, \beta_{-2})\\
    		A_{4}(\varphi,\mathbf{n}_x, \beta_{-2}, \Gamma_{\mathbf{n},3,\phi})\\
    		A_{3}(\varphi, \Gamma_{\mathbf{n},4,\phi}, \Gamma_{\mathbf{n},3,\phi}, \Gamma_{\mathbf{n},3,\phi})\\
    		A_{3}(\varphi, \Gamma_{\mathbf{n},4,\phi}, \Gamma_{\mathbf{n},3,\phi}, \beta_{-1})\\
    		A_{3}(\varphi, \Gamma_{\mathbf{n},4,\phi}, \Gamma_{\mathbf{n},3,\phi}, \beta_{0})\\
    	\end{pmatrix}
\end{aligned}
\end{equation}
% \begin{equation}
%     \begin{pmatrix}
%         f_{-4} \\
%         f_{2} \\
%         f_{31}
%     \end{pmatrix} = \begin{pmatrix}
%         1 & 0 & 0 \\
%         - 1 & 1 & 0 \\
%         0 & - 1 & 1
%     \end{pmatrix}\begin{pmatrix}
%     		A_{3\mathrm{rel}}(\varphi,\betaGo,\betaGtwo,\beta_{1}) \\
%     		A_{3\mathrm{rel}}(\varphi,\betaGo,\betaGtwo,\beta_{2}) \\
%     		A_{3\mathrm{rel}}(\varphi,\betaGo,\betaGtwo,\betaGtwo)
%     	\end{pmatrix},
% \end{equation}
% \begin{equation}
%     \begin{pmatrix}
%         f_{32} \\
%         f_{4} \\
%         f_{51}
%     \end{pmatrix} = \begin{pmatrix}
%         1 & 0 & 0 \\
%         0 & 1 & 0 \\
%         0 & 0 & 1
%     \end{pmatrix}\begin{pmatrix}
%         A_{4\mathrm{rel}}\left( \varphi,n_{x},\betaGtwo,\beta_{3} \right) \\
%         A_{4\mathrm{rel}}\left( \varphi,n_{x},\beta_{3},\beta_{4} \right) \\
%         A_{4\mathrm{rel}}\left( \varphi,n_{x},\beta_{4},\betaGth \right)
%     \end{pmatrix},
% \end{equation}
% and
% \begin{equation}
% 	\begin{pmatrix}
% 		f_{52} \\
% 		f_{6}
% 		\end{pmatrix} = \begin{pmatrix}
% 		1 & - 1 \\
% 		0 & 1
% 		\end{pmatrix}\begin{pmatrix}
% 		A_{3\mathrm{rel}}(\varphi,\betaGf,\betaGth,\betaGth) \\
% 		A_{3\mathrm{rel}}(\varphi,\betaGf,\betaGth,\beta_{5})
% 	\end{pmatrix}
% \end{equation}
with \( w_{-3,\mathbf{n}}=w_{-3,\mathbf{n}}^1+w_{-3,\mathbf{n}}^2\) and analogously \( w_{-2,\mathbf{n}} = w_{-2,\mathbf{n}}^1 + w_{-2,\mathbf{n}}^2\).
% This example clearly demonstrates the general algorithm for generating the two-dimensional weighting factors.
\begin{algorithm}[t]%
    \caption{Computation of \(\SystemMatrix\) in 2D.}%
    \label{algo:forward projection 2d}
    \DontPrintSemicolon
    \KwData{Image \( \mathbf{u} \) to be projected}
    Initialize $\mathbf{p} = \mathbf{0}\in \R^{N_p\times N_d}$.\;
    \For{projection angles \(\phi\)}{
    \For{voxels \(\mathbf{n}\)}{
        Compute angles \((\gamma_{\mathbf{n},i,\phi})_{i=1}^4\) with \eqref{eq:gamma_computation}.\;
        Compute relevant detectors \( \mathcal{D} \) with \cref{eq:relevant detectors 2d}.\;
        Compute angles \( (\beta_\mathbf{m})_{\mathbf{m}\in \mathcal{D}} \) with \eqref{eq:beta_computation}.\;
        Compute weight $w_{\mathbf{m},\mathbf{n}}$ according to the procedure \eqref{eq:2d procedure}.\;
        \For{\(\mathbf{m}\in\mathcal{D}\)}{
            \(\mathbf{p}_{\phi,\mathrm{m}} \mathrel{+}= w_{\mathbf{m},\mathbf{n}} \mathbf{u}_{\mathbf{n}}\).
        }
        }
    }
    \KwResult{Projections \( \mathbf{p} = \SystemMatrix \mathbf{u} \)}
\end{algorithm}
% \begin{algorithm}[t]%
%     \caption{TODO Computation of the backprojection in two dimensions.}%
%     \label{algo:backprojection 2d}
%     \DontPrintSemicolon
%     \KwData{Image \( \mathbf{u} \) to be projected}
%     \For{pixels \( k = 1, 2, \dotsc, N\)}{
%         Compute relevant detectors \( \mathcal{D} \) with \cref{eq:relevant detectors 2d}.\;
%         \For{relevant detectors \( p \in \mathcal{D} \)}{
%             Compute angles \( \beta \) ... and depending on their relation compute the area factors .... and set \( \mathbf{p}_p = \) ...
%         }
%     }
%     \KwResult{Projections \( \mathbf{p} = \SystemMatrix \mathbf{u} \)}
% \end{algorithm}
\subsection{The Three-dimensional Case}%
\label{three-dimensional-case}
Since we assume the rotation $\phi$ only around the $z$ axis, from the top view the shapes emerging in the 2D case stay largely the same. Moreover, we can still use the angles $\gamma_{\mathbf{n},i,\phi}$ and $\beta_{\mathbf{m}_y}$ as before.
However, it is now necessary to determine in addition to the \( y \) component the \( z \) position of the involved detector elements.
From~\cref{eq:z expr 1}, it can be derived
\begin{equation}
\begin{split}
	\frac{s + d}{d_{z}}\frac{\cos{(\varphi} - \gamma_\mathrm{min})}{\cos\gamma_\mathrm{min}}&\frac{z_{\min}}{s\cos\varphi - x} < \mathbf{m}_{z}\\ 
    < &\frac{s+d}{d_{z}}\frac{\cos{(\varphi} - \gamma_\mathrm{max})}{\cos\gamma_\mathrm{max}}\frac{z_{\max}}{s\cos\varphi - x}.
\end{split}
\end{equation}
where \( z_{max} \) and \( z_{min} \) are chosen to be the upper and lower limits of the voxel element, respectively and \(\gamma_\mathrm{min}\), \(\gamma_\mathrm{max}\) the corresponding angles from the top view.
In the case of \( V_{\mathbf{n},02} \) and \( V_{\mathbf{n},12} \), \eg, the following relations are valid, considering~\cref{eq:cos phi nx a}
\begin{equation}
\begin{split}
	\frac{s+d}{d_{z}}&\frac{\sin{(\gamma_{\mathbf{n},\phi,1}} - \gamma_{\mathbf{n},\phi,2})}{\sin{(\varphi} - \gamma_{\mathbf{n},\phi,1})\cos\gamma_{\mathbf{n},\phi,2}}\frac{\mathbf{n}_{z}c}{a} < \mathbf{m}_{z} \\
    < &\frac{s+d}{d_{z}}\frac{\sin{(\gamma_{\mathbf{n},\phi,1}} - \gamma_{\mathbf{n},\phi,2})}{\sin{(\varphi} - \gamma_{\mathbf{n},\phi,1})\cos\gamma_{\mathbf{n},\phi,2}}\frac{\left( \mathbf{n}_{z} + 1 \right)c}{a}.
\end{split}
\end{equation}
In the following we will derive explicit formulae for the possible sub-voxel volumes emerging from intersecting the rays to a detector element with a voxel. 
% As this is usually satisfied in practice, we will assume that in $z$-direction each voxel contributes to at most three rows of detector elements. 
% Let us denote three vertically stacked sub-voxel volumes corresponding to three likewise vertically stacked detector elements as $\vol_\mathrm{top}$, $\vol_\mathrm{bottom}$, and $\vol_\mathrm{center}$. It is sufficient to derive a formula for $\vol_\mathrm{top}$, since $\vol_\mathrm{bottom}$ can be treated analogously and $\vol_\mathrm{center} = (\vol_\mathrm{top}+\vol_\mathrm{center}+\vol_\mathrm{bottom})-(\vol_\mathrm{top}-\vol_\mathrm{bottom})$ where the first parenthesis is a simple prism, \cf~\cref{fig:3d_composition}.

In the 3D setting the case distinction to be made is more involved than in 2D. As the comprehensive derivation of all formulae would be beyond the scope of this manuscript \footnote{and, in particular, beyond the page-limit}, we restrict to several instructive cases and refer the interested reader to the published source code\footnote{Interested readers may of course also contact us directly.}.
% The following examples are restricted to the computation of the top sub-volume (\ie, with largest $z$ component) of a voxel.
% \textcolor{red}{Andi: Ich bin unsicher ob meine Beschriebung der case distinction passt, aber wir müssen irgendeine guideline für den Leser geben. Es ist sonst zu schwierig, systematisch zu verstehen, wie man alle subvoxel abdeckt kann (zumindest theoretisch) und warum gerade die aufgelisteten gezeigt werden...}
\subsubsection{Top-view trapezoid}
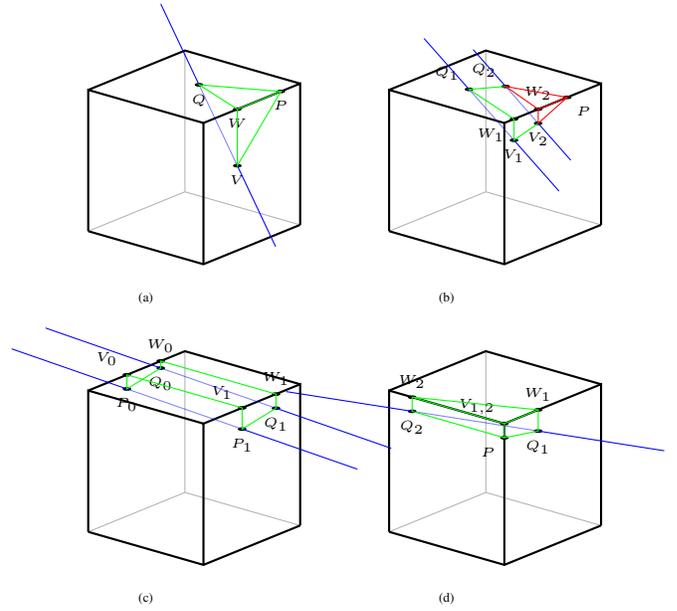
\begin{figure}
	\centering
    \begin{tikzpicture}[font=\tiny]
        \begin{scope}[3d view={40}{20}]
            \coordinate (a) at (0,-1,-1.5);
            \node at (a) [below] {(a)};
            \drawcubeTD{-1}{-1}{-1}{2};
            \coordinate (V) at (1, -0.3, 0.2);
            \coordinate (Q) at (0, 0.1, 1);
            \draw[fill] (Q) circle (0.05);
            \draw[fill] (V) circle (0.05);
            \node at (Q) [below] {\( Q \)};
            \node at (V) [below] {\( V \)};
            \opaqueconnect[blue]{Q}{V};
            \coordinate (P) at (1,0.6,1);
            \node at (P) [below] {\( P \)};
            \coordinate (W) at (1,-0.3,1);
            \node at (W) [below] {\( W \)};
            \draw[fill] (P) circle (0.05);
            \draw[fill] (W) circle (0.05);
            \draw [green] (W) -- (P) -- (Q) -- (W) --(V) -- (P);
        \end{scope}
        \begin{scope}[xshift=4cm,3d view={40}{20}]
            \coordinate (b) at (0,-1,-1.5);
            \node at (b) [below] {(b)};
            \drawcubeTD{-1}{-1}{-1}{2};
            \coordinate (W1) at (1, -0.8, 1);
            \node at (W1) [below left] {\( W_1 \)};
            \draw[fill] (W1) circle (0.05);
            \coordinate (W2) at (1, -0.3, 1);
            \node at (W2) [above] {\( W_2 \)};
            \draw[fill] (W2) circle (0.05);
            \coordinate (P) at (1, 0.3, 1);
            \node at (P) [below right] {\( P \)};
            \draw[fill] (P) circle (0.05);
            \coordinate (Q1) at (-0.2, -0.3, 1);
            \node at (Q1) [above left] {\( Q_1 \)};
            \draw[fill] (Q1) circle (0.05);
            \coordinate (Q2) at (0.1, 0.1, 1);
            \node at (Q2) [above left] {\( Q_2 \)};
            \draw[fill] (Q2) circle (0.05);
            \coordinate (V1) at (1, -0.8, 0.7);
            \coordinate (V2) at (1, -0.3, 0.8);
            \node at (V1) [below] {\( V_1 \)};
            \draw[fill] (V1) circle (0.05);
            \node at (V2) [below] {\( V_2 \)};
            \draw[fill] (V2) circle (0.05);
            % \node at (V1) [left]{$V_1$};
            % \node at (V2) [right]{$V_2$};
            % \draw[fill] (V1) circle (0.05);
            % \draw[fill] (V2) circle (0.05);
            \opaqueconnect[blue]{Q1}{V1};
            \opaqueconnect[blue]{Q2}{V2};
            \draw [green] (W1) -- (Q1) -- (Q2);
            \draw [green] (W1) -- (V1) -- (V2);
            \draw [red] (P) -- (Q2) -- (W2) -- (V2) -- (P) -- (W2);
        \end{scope}
        \begin{scope}[xshift=4cm, yshift=-4cm,3d view={40}{20}]
            \coordinate (d) at (0,-1,-1.5);
            \node at (d) [below] {(d)};
            \drawcubeTD{-1}{-1}{-1}{2};
            \coordinate (W1) at (1, -0.3, 1);
            \node at (W1) [above]{$W_1$};
            \draw[fill] (W1) circle (0.05);
            \coordinate (Q1) at (1, -0.3, 0.7);
            \node at (Q1) [below]{$Q_1$};
            \draw[fill] (Q1) circle (0.05);
            \coordinate (W2) at (-0.6, -1, 1);
            \node at (W2) [above]{$W_2$};
            \draw[fill] (W2) circle (0.05);
            \coordinate (Q2) at (-0.6, -1, 0.8);
            \node at (Q2) [below]{$Q_2$};
            \draw[fill] (Q2) circle (0.05);
            \coordinate (P) at (1,-1, 0.8);
            \node at (P) [below left]{$P$};
            \draw[fill] (P) circle (0.05);
            \opaqueconnect[blue]{Q1}{Q2};
            \coordinate(V12) at (1,-1,1);
            \node at (V12) [above left]{$V_{1,2}$};
            \draw[fill] (V12) circle (0.05);
            \draw [green] (P) -- (Q2) --(W2)--(W1)--(Q1)--(P)--(V12)--(W2);
            \draw [green] (V12)--(W1)--cycle;
        \end{scope}
        \begin{scope}[xshift=0cm, yshift=-4cm,3d view={40}{20}]
            \coordinate (c) at (0,-1,-1.5);
            \node at (c) [below] {(c)};
            \drawcubeTD{-1}{-1}{-1}{2};
            \coordinate (W0) at (-1, 0.5, 1);
            \node at (W0) [above] {\( W_0 \)};
            \draw[fill] (W0) circle (0.05);
            \coordinate (W1) at (1, 0.5, 1);
            \node at (W1) [above] {\( W_1 \)};
            \draw[fill] (W1) circle (0.05);
            \coordinate (V0) at (-1, -0.2, 1);
            \node at (V0) [above left] {\( V_0 \)};
            \draw[fill] (V0) circle (0.05);
            \coordinate (V1) at (1, -0.2, 1);
            \node at (V1) [above left] {\( V_1 \)};
            \draw[fill] (V1) circle (0.05);
            \coordinate (P0) at (-1, -0.2, 0.8);
            \node at (P0) [below] {\( P_0 \)};
            \draw[fill] (P0) circle (0.05);
            \coordinate (P1) at (1, -0.2, 0.7);
            \node at (P1) [below] {\( P_1 \)};
            \draw[fill] (P1) circle (0.05);
            \coordinate (Q0) at (-1, 0.5, 0.9);
            \node at (Q0) [below] {\( Q_0 \)};
            \draw[fill] (Q0) circle (0.05);
            \coordinate (Q1) at (1, 0.5, 0.8);
            \node at (Q1) [below] {\( Q_1 \)};
            \draw[fill] (Q1) circle (0.05);
            \opaqueconnect[blue]{P0}{P1};
            \opaqueconnect[blue]{Q0}{Q1};
            \draw [green] (V0) -- (V1) --(P1)--(Q1)--(W1)--(W0)--(Q0)--(P0)--cycle;
        \end{scope}
    \end{tikzpicture}
	\caption{
		Various cases of contribution of a sub-voxel (delimited by green lines) to a detector element.
		Rays to corners of the detector element are shown in blue.
		(a) Sub-voxel VPQW with the coefficient vector \( \mathbf{C_{m}} = (1, 0, 0, 0) \).
		(b) Sub-voxel \( V_{1}V_{2}W_{2}Q_{2}Q_{1}W_{1} \) with the coefficient vector \( \mathbf{C_{m}} = (1, 1, 0, 0) \).
		The volume is the difference between the volumes \( V_{1}PQ_{1}W_{1} \) and \( V_{2}PW_{2}Q_{2} \).
		(c) Sub-voxel \( P_{1}Q_{1}Q_{0}P_{0} \) \( V_{1}W_{1}V_{0}W_{0} \) with the coefficient vector \( \mathbf{C_{m}} = (1, 1, 1, 1) \).
		(d) Sub-voxel \( PQ_{1}Q_{2}V_{12}W_{1}W_{2} \) with the coefficient vector \( \mathbf{C_{m}} = (1, 1, 1) \).}
    \label{fig:3d_cases}
\end{figure}

We first cover cases, where from a top perspective the intersection of the area contained within the two rays delimiting the sub-voxel in $y$ direction and the top surface of the sub-voxel comprises a trapezoid. In the following we distinguish sub-cases depending on wether each of the rays intersects the front and/or rear surface of the voxel\footnote{front/rear from the perspective of the source}.
\paragraph{Tetrahedron}
Consider the sub-case that only one of the two delimiting rays intersects the voxel and only either at the front or the rear surface(\cf~\cref{fig:3d_cases} (a)) leading to a tetrahedron, 
The resulting volume \( \vol \) can be calculated by determining the points, where rays to the upper boundary of the detector hit its top surface.

With \( V = ((\mathbf{n}_{x} + 1)a, y, z)^T \) and using~\cref{eq:x expr 1} and~\cref{eq:y expr 1}, the vectors to these points \( \mathbf{p} \) and \( \mathbf{q} \) can be written as
\begin{equation}
\begin{aligned}
	\mathbf{p} = \begin{pmatrix}
		\left( \mathbf{n}_{x} + 1 \right)a \\
		s\sin\varphi + \frac{1}{\sin\varphi} \ \bigl\lbrack \left( s\cos\varphi - \left( \mathbf{n}_{x} + 1 \right)a \right)\cos\varphi
        \\ -\frac{s+d}{\mathbf{m}_{z}d_{z}}\left( \mathbf{n}_{z} + 1 \right)c\bigr\rbrack \\
		\left( \mathbf{n}_{z} + 1 \right)c
	\end{pmatrix}
 \label{eq:pyramid p}
 \end{aligned}
\end{equation}
and
\begin{equation}
	\mathbf{q} = \begin{pmatrix}
		s\cos\varphi - \frac{s + d}{\mathbf{m}_{z}d_{z}}\frac{\cos{(\varphi} - \beta_{\mathbf{m}_y})}{\cos\beta_{\mathbf{m}_y}}\left( \mathbf{n}_{z} + 1 \right)c \\
		s\sin\varphi - \frac{s + d}{\mathbf{m}_{z}d_{z}}\frac{\sin{(\varphi} - \beta_{\mathbf{m}_y})}{\cos\beta_{\mathbf{m}_y}}\left( \mathbf{n}_{z} + 1 \right)c \\
		\left( \mathbf{n}_{z} + 1 \right)c
	\end{pmatrix}.
  \label{eq:pyramid q}
\end{equation}
With~\cref{eq:z expr 1} and~\cref{eq:y expression}, the vector \( \mathbf{v} \) to the point \( V \) can be expressed as
\begin{equation}
	\mathbf{v} = \begin{pmatrix}
		\left( \mathbf{n}_{x} + 1 \right)a \\
		s\sin\varphi - (s\cos\varphi - \left( \mathbf{n}_{x} + 1 \right)a)\tan{(\varphi} - \beta_{\mathbf{m}_y})\  \\
		\frac{\mathbf{m}_{z}d_{z}}{s + d}\frac{\cos\beta_{\mathbf{m}_y}}{\cos{(\varphi} - \beta_{\mathbf{m}_y})}\left( s\cos\varphi - \left( \mathbf{n}_{x} + 1 \right)a \right)
	\end{pmatrix}.
   \label{eq:pyramid v}
\end{equation}
The volume of the resulting tetrahedron can be calculated as
\begin{equation}
	{\vol} = \frac{1}{6}\left| \left( \left( \mathbf{n}_{x} + 1 \right)a - \mathbf{q}_{x} \right)\left( \mathbf{v}_{y} - \mathbf{p}_{y} \right)\left( \left( \mathbf{n}_{z} + 1 \right)c - \mathbf{v}_{z} \right) \right|.
\end{equation}
Using~\cref{eq:pyramid p} to~\cref{eq:pyramid v}, this can be written as
\begin{equation}
	{\vol} = \frac{1}{6}a^{3}\left| \frac{\tan\alpha_{\mathbf{m}_{z}}}{\sin\varphi}\ \frac{\cos{(\varphi} - \beta_{\mathbf{m}_y})}{\cos\beta_{\mathbf{m}_y}} g_{\mathbf{m}_y}^{3} \right|,
\end{equation}
where
\begin{equation}
\begin{split}
	&g_{\mathbf{m}_y} = \frac{c}{a}\frac{\mathbf{n}_{z} + 1}{\tan\alpha_{\mathbf{m}_{z}}}- \frac{\cos{\beta_{\mathbf{m}_y}}}{\cos{(\varphi-\beta_{\mathbf{m}_y})}}\frac{s\cos\varphi - \left( \mathbf{n}_{x} + 1 \right)a}{a} \\ 
    &=\frac{c}{a}\frac{\mathbf{n}_{z} + 1}{\tan\alpha_{\mathbf{m}_{z}}}- \frac{\cos{\beta_{\mathbf{m}_y}}\cos{(\varphi} - \gamma_{\mathbf{n},2,\phi})\sin{(\varphi} - \gamma_{\mathbf{n},1,\phi})}{\cos{(\varphi-\beta_{\mathbf{m}_y})}\sin(\gamma_{\mathbf{n},1,\phi} - \gamma_{\mathbf{n},2,\phi})}
\end{split}
\end{equation}
and \(\tan\alpha_{\mathbf{m}_{z}} = \frac{\mathbf{m}_{z}d_{z}}{s + d}\).
Analogously, if the the intersecting point is \( V = (\mathbf{n}_{x}a, y, z)^T \) the resulting volume \( \vol \) can be calculated as
\begin{equation}
	{\vol} = \frac{1}{6}a^{3}\left| \frac{\tan\alpha_{\mathbf{m}_{z}}}{\sin\varphi}\ \frac{\cos{(\varphi} - \beta_{\mathbf{m}_y})}{\cos\beta_{\mathbf{m}_y}} f_{\mathbf{m}_y}^{3} \right|,
\end{equation}
where
\begin{equation}
\begin{split}
	&f_{\mathbf{m}_y} = \frac{c}{a}\frac{\mathbf{n}_{z} + 1}{\tan\alpha_{\mathbf{m}_{z}}}- \frac{\cos{\beta_{\mathbf{m}_y}}}{\cos{(\varphi-\beta_{\mathbf{m}_y})}}\frac{s\cos\varphi - \mathbf{n}_{x}a}{a} \\ 
    &= \frac{c}{a}\frac{\mathbf{n}_{z} + 1}{\tan\alpha_{\mathbf{m}_{z}}} - \frac{\cos{\beta_{\mathbf{m}_y}}\cos{(\varphi} - \gamma_{\mathbf{n},1,\phi})\sin{(\varphi} - \gamma_{\mathbf{n},2,\phi})}{\cos{(\varphi-\beta_{\mathbf{m}_y})}\sin(\gamma_{\mathbf{n},1,\phi} - \gamma_{\mathbf{n},2,\phi})}.
\end{split}
\end{equation}
For the factors \( f_{\mathbf{m}_y} \) and \( g_{\mathbf{m}_y} \), the following relation is valid,
\begin{equation}
	g_{\mathbf{m}_y} - f_{\mathbf{m}_y} = \frac{\cos{\beta_{\mathbf{m}_y}}}{\cos{(\varphi-\beta_{\mathbf{m}_y})}}.
 \label{eq:f minus g}
\end{equation}

\paragraph{Difference of Tetrahedra}
Next consider the sub-case that both rays intersect the voxel, albeit only at the front or at the rear surface of the voxel, \cf~\cref{fig:3d_cases}, (b).
Let us denote \( V_{1} = (({\mathbf{n}_{x}} + 1){a}, {y_{1}}, {z_{1}})^{T} \) and \( V_{2} = (({\mathbf{n}_{x}} + 1){a}, {y_{2}}, {z_{2}})^{T} \). 
The resulting volume can be regarded as a composition of the individual tetrahedral volumes of both rays.
It can be obtained by subtracting the smaller from the larger volume leading to
\begin{equation}
\begin{split}
	\vol = \frac{1}{6}a^{3}\left| \frac{\tan\alpha_{\mathbf{m}_{z}}}{\sin\varphi}\left\lbrack \frac{\cos{(\varphi-\beta_{\mathbf{m}})}}{\cos{\beta_{\mathbf{m}}}}g_{\mathbf{m}}^{3} - \right.\right.
    \\ \left.\left. -\frac{\cos{(\varphi-\beta_{\mathbf{m}+1})}}{\cos{\beta_{\mathbf{m}+1}}}g_{\mathbf{m} + 1}^{3} \right\rbrack \right|.
\end{split}
\end{equation}
\paragraph{Prismatoid}
If the sub-voxel is determined by two rays which do not intersect to top surface as in~\cref{fig:3d_cases} (c), the resulting volume of the prismatoid is computed as

% \begin{equation}
% 	\begin{aligned}
% 		\vol = \frac{a}{6}\left\lbrack 2\frac{\left| V_{1}P_{1} \right| + \left| W_{1}Q_{1} \right|}{2}\left| V_{1}W_{1} \right| + \right. \\
%         + \frac{\left| V_{1}P_{1} \right| + \left| W_{1}Q_{1} \right|}{2}\left| V_{0}W_{0} \right| +\\
% 		+ \frac{\left| V_{0}P_{0} \right| + \left| W_{0}Q_{0} \right|}{2}\left| V_{1}W_{1} \right| + \\
%         \left. + 2\frac{\left| V_{0}P_{0} \right| + \left| W_{0}Q_{0} \right|}{2}\left| V_{0}W_{0} \right| \right\rbrack.
% 	\end{aligned}
% \end{equation}
% In further consequence
% \begin{equation}
% \begin{split}
% 	\vol = \frac{1}{6}a^{3}\ \left| \frac{\tan\alpha_{\mathbf{m}_{z}}}{\sin\varphi}\lbrack \left( f_{\mathbf{m}_y}^{2} + f_{\mathbf{m}_y}g_{\mathbf{m}_y} + g_{\mathbf{m}_y}^{2} \right) - \right.\\
%     \left.- \left( f_{\mathbf{m}_y + 1}^{2} + f_{\mathbf{m}_y + 1}g_{\mathbf{m}_y + 1} + g_{\mathbf{m}_y + 1}^{2} \right) \rbrack \ \right|.
%  \end{split}
% \end{equation}
% Using the relation~\cref{eq:f minus g}, the final result is:
\begin{equation}
\begin{split}
	\vol = &\frac{1}{6}a^{3}\bigg\lvert \frac{\tan\alpha_{\mathbf{m}_{z}}}{\sin\varphi}\left\lbrack \frac{\cos{(\varphi-\beta_{\mathbf{m}_y})}}{\cos{\beta_{\mathbf{m}_y}}}\left( g_{\mathbf{m}_y}^{3} - f_{\mathbf{m}_y}^{3} \right) - \right. \\
    &\left. - \frac{\cos{(\varphi-\beta_{\mathbf{m}_y+1})}}{\cos{\beta_{\mathbf{m}_y+1}}}{\left( g_{\mathbf{m}_y + 1}^{3} - f_{\mathbf{m}_y + 1}^{3} \right)}\right\rbrack \bigg\rvert.
\end{split}
\end{equation}
Astonishingly, it turns out that we can derive a general formula for the relative volume (in relation to the orthorhombic voxel volume \( \vol_{\mathrm{voxel}} = abc) \) \( \vol_{\mathrm{rel}} \) in the case \emph{top-view trapezoid} which reads as
\begin{equation}
\begin{split}
	&\vol_\mathrm{rel}(\mathbf{m}_{z},\beta_{\mathbf{m}_y},\ \beta_{\mathbf{m}_y + 1},\mathbf{C}_{\mathbf{m}_{\mathbf{z}}})\\
    = &\frac{1}{6}\frac{a^{2}}{bc}\Biggl | \frac{\tan\alpha_{\mathbf{m}_{z}}}{\sin\varphi} \Biggl \lbrack \frac{\cos{(\varphi-\beta_{\mathbf{m}_y})}}{\cos{\beta_{\mathbf{m}_y}}}\left( C_{g,\mathbf{m}_y}g_{\mathbf{m}_y}^{3} - C_{f,\mathbf{m}_y}f_{\mathbf{m}_y}^{3} \right) \\ 
    &- {\frac{\cos{(\varphi-\beta_{\mathbf{m}_y+1})}}{\cos{\beta_{\mathbf{m}_y+1}}}\left( C_{g,\mathbf{m}_y+1}g_{\mathbf{m}_y + 1}^{3} - C_{f,\mathbf{m}_y+1}f_{\mathbf{m}_y + 1}^{3} \right)}\Biggr \rbrack \Biggr |,
 \end{split}
 \label{eq:rel vol 3d}
\end{equation}
where we define \(\mathbf{C}_{\mathbf{m}} = (C_{g,\mathbf{m}_y}, C_{f,\mathbf{m}_y}, C_{g,\mathbf{m}_y+1}, C_{f,\mathbf{m}_y+1})\) with \(C_{g,\mathbf{m}_y} = 1\) if the ray to detector element $\mathbf{m}$ has an intersection with the front surface of the voxel and zero else and analogously for \(C_{f,\mathbf{m}_y}\) and the rear surface.
The values zero or one are intrinsically related to the signs of the associated functions \( g_{\mathbf{m}_y} \) and \( f_{\mathbf{m}_y} \).
% For convenience, in the function \( \vol_\mathrm{rel} \) only the variables that change within a single voxel are listed.
The resulting volume elements can be calculated directly or by the negative combination of two volume elements
\begin{equation}
\begin{aligned}
	\vol_{\mathrm{res}}(\mathbf{m}_{z},\beta_{\mathbf{m}_y},\ &\beta_{\mathbf{m}_y + 1}) \\
    =&\vol_{\mathrm{rel}}(\mathbf{m}_{z} - 1,\beta_{\mathbf{m}_y},\ \beta_{\mathbf{m}_y + 1},\mathbf{C}_{\mathbf{m}^{-1}}) \nonumber \\
    &-\vol_{\mathrm{rel}}(\mathbf{m}_{z},\beta_{\mathbf{m}_y},\ \beta_{\mathbf{m}_y + 1}\mathbf{C}_{\mathbf{m}}).
\end{aligned}
\end{equation}
with \(C_{\mathbf{m^{-1}}}=C_{(\mathbf{m}_y,\mathbf{m}_z-1)}\) the vector for the sub-volume corresponding the the detector element at $(\mathbf{m}_y,\mathbf{m}_z-1)$.
At the singularity \( \varphi = 0 \), considering that \( g_{\mathbf{m}_y+1} = g_{\mathbf{m}_y} \), \( f_{\mathbf{m}_y+1} = f_{\mathbf{m}_y} \), \( C_{g,\mathbf{m}_y+1} = C_{g,\mathbf{m}_y} \), and \( C_{f,\mathbf{m}_y+1} = C_{f,\mathbf{m}_y} \),~\cref{eq:rel vol 3d}  can be expressed as:
\begin{equation}
\begin{split}
	&\vol_\mathrm{rel}(\mathbf{m}_{z},\beta_{\mathbf{m}_y},\ \beta_{\mathbf{m}_y + 1},\mathbf{C}_{\mathbf{m}})
    \\ = &\frac{1}{6}\frac{a^{2}}{bc}\Biggl | \tan\alpha_{\mathbf{m}_{z}}\Biggl \lbrack C_{g,\mathbf{m}_y}g_{\mathbf{m}_y}^{2}\left( g_{\mathbf{m}_y} + 3\frac{s}{a} - 3\mathbf{n}_{x} - 3 \right)\ \\
    &-C_{f,\mathbf{m}_y}f_{\mathbf{m}_y}^{2}\left( f_{\mathbf{m}_y} + 3\frac{s}{a} - 3\mathbf{n}_{x} \right)\ \Biggr \rbrack\left(\tan\beta_{\mathbf{m}_y} - \tan\beta_{\mathbf{m}_y + 1} \right)\Biggr |.
 \end{split}
 \label{eq:vol}
 \end{equation}
 
\subsubsection{Top-view triangle}
Secondly let us consider the case that from the top perspective the intersection of the area between the two delimiting rays and the voxel top surface is a triangle delimitted only by one ray and the voxel boundaries, \cf~\cref{fig:3d_cases} (d). A similar strategy as in the previous case (\cf~\cref{sec:appendix chopped prism}) leads to
% \begin{equation}
% \begin{split}
% 	\vol = &a^{3}\left| \frac{\tan\alpha_{m_{z}}}{\sin\varphi}\left\lbrack\frac{\cos{(\varphi-\beta_{i})}}{\cos{\beta_{i}}} g_{i}^{3} + \right. \right.
%     \\ &\left.\left. + \tan\varphi\frac{\sin{(\varphi-\beta_{i})}}{\cos{\beta_{i}}} h_{i}^{3}- \frac{1}{cos\varphi}g_{g}^{3}\right\rbrack \right|.
% \end{split}
% \end{equation}
% Other possible polyhedrons are treated in the supplementary information. In general, the relative volume (in relation to the orthorhombic voxel
% volume \( \vol_{\mathrm{rel}} = abc \)) \( \vol_{\mathrm{voxel}} \) can be expressed as:
\begin{equation}
\begin{split}
	&\vol_{\mathrm{rel}}(\mathbf{m}_{z},\beta_{\mathbf{m}_y},\mathbf{C}_{\mathbf{m}})\\
    =&\frac{1}{6}\frac{a^{2}}{bc}\Biggl | \frac{\tan\alpha_{\mathbf{m}_{z}}}{\sin\varphi}\Biggl \lbrack C_{y,\mathbf{m}_y}\frac{\cos{(\varphi-\beta_{\mathbf{m}_y})}}{\cos{\beta_{\mathbf{m}_y}}}g_{\mathbf{m}_y}^{3}\\
    &+  C_{x,\mathbf{m}_y}\tan\varphi\frac{\sin{(\varphi-\beta_{\mathbf{m}_y})}}{\cos{\beta_{\mathbf{m}_y}}} h_{\mathbf{m}_y}^{3} - C_{g}\frac{1}{cos\varphi}g_{g}^{3}\Biggr \rbrack \Biggr |,
\end{split}
\label{eq:res vol 3d}
\end{equation}
where, similar to the previous elaborations, the components of the vector \( \mathbf{C}_{\mathbf{m}} = (C_{y,\mathbf{m}_y}, C_{x,\mathbf{m}_y}, C_{g})^T \) are either 1 or 0, if intersection occurs or not, respectively and we refer to the appendix for details.
% As in~\cref{eq:vol}, the values zero/one are related to the signs of the associated functions \( g_\mathbf{m}_y \), \( g_{iy} \) and \( g_g \).
The resulting volume elements can be calculated directly or by the negative combination of two volume elements:
\begin{equation}
\begin{split}
	&\vol_\mathrm{res}(\mathbf{m}_{z},\beta_{\mathbf{m}_y},\mathbf{C}_{\mathbf{m}})\\ 
    &= \vol_\mathrm{rel}(\mathbf{m}_{z} - 1,\beta_{\mathbf{m}_y},\ \mathbf{C}_{\mathbf{m}^{-1}}) - \vol_\mathrm{rel}(\mathbf{m}_{z},\beta_{\mathbf{m}_y},\mathbf{C}_{\mathbf{m}}).
\end{split}  
\end{equation}
At the singularity \( \varphi = 0 \), the above can be expressed as:
\begin{equation}
\begin{split}
	&\vol_\mathrm{rel}(\mathbf{m}_{z},\beta_{\mathbf{m}_y},\ \mathbf{C}_{\mathbf{m}}) \\
    = &\frac{1}{6}\frac{a^{2}}{bc}\Biggl |\tan\alpha_{\mathbf{m}_{z}}\Biggl \lbrack C_{y,\mathbf{m}_y}g_{\mathbf{m}_y}^{2}\left(\tan{\beta_{\mathbf{m}_y}} g_{\mathbf{m}_y} \right.\\
    &\left.- 3\frac{\sin\gamma_{\mathbf{n},1,\phi}\sin{(\beta_{\mathbf{m}_y}-\gamma_{\mathbf{n},2,\phi}})}{\sin{(\gamma_{\mathbf{n},1,\phi}} - \gamma_{\mathbf{n},2,\phi})\cos\beta_{\mathbf{m}_y}}\right)- C_{x,\mathbf{m}_y}\tan\beta_{\mathbf{m}_y}h_{\mathbf{m}_y}^{3} \Biggr \rbrack \Biggr|.
 \end{split}
\end{equation}
Note that \(g_{g} = g_{\mathbf{m}_y}\) (see ~\cref{eq:gi gg} ) and  \(C_{g} = C_{y,\mathbf{m}_y}\).
\begin{figure}
	\centering
	\includegraphics[width=0.25\textwidth]{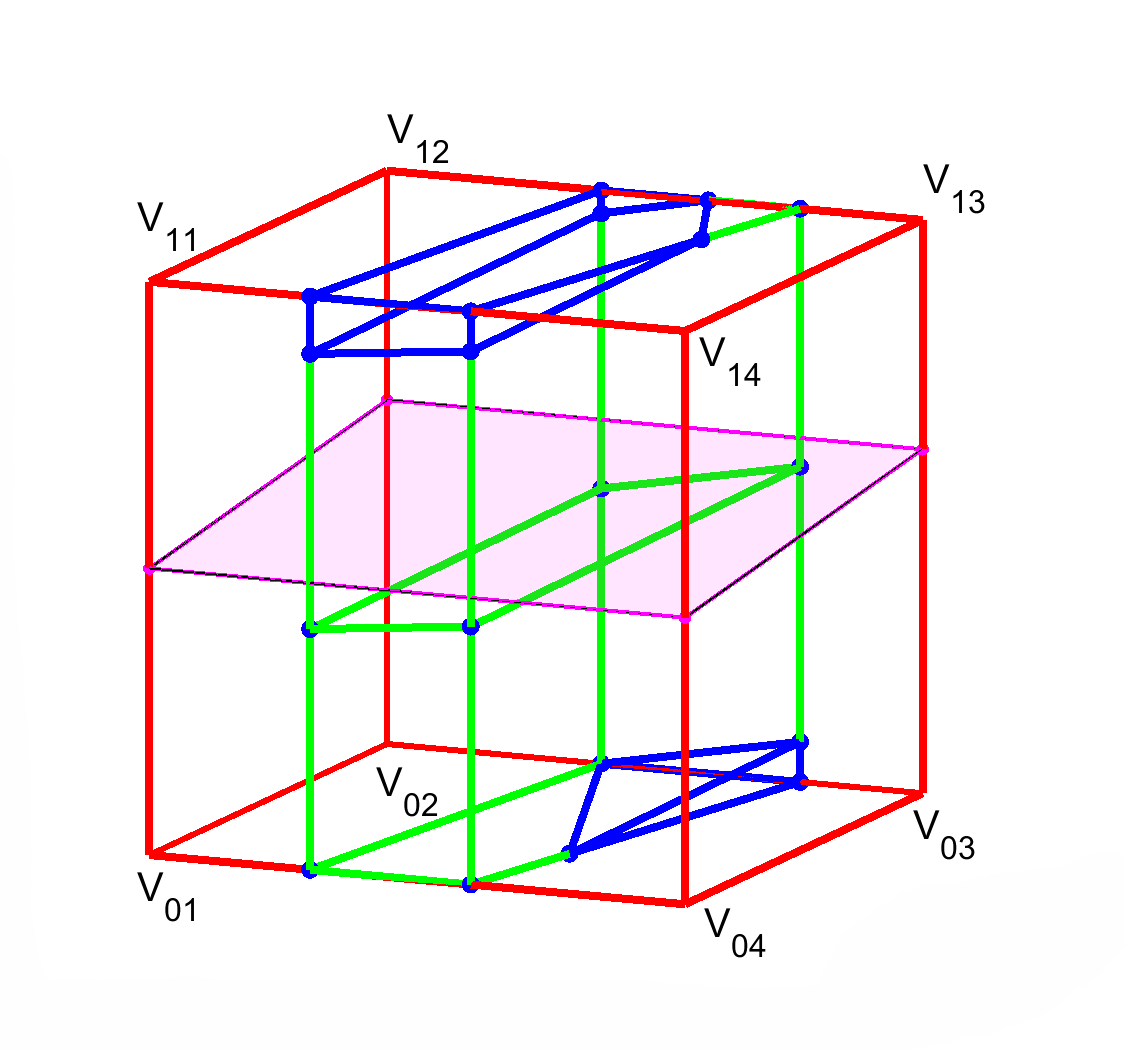}
	\caption{
		Determining the volume of sub-voxels contributing to a detector element by proper segmentation.
		The volume of the blue delimited sub-voxels with the coefficient vectors \( \mathbf{C_{m+1}}=(1,1,0,1) \) and \( \mathbf{C_{m}}=(0,0,1,0) \) can be calculated directly.
		The volume of the green delimited sub-voxels are calculated by subtraction.
		The plane in magenta color splits the voxel in an upper and lower part.
	}%
    \label{fig:3d_composition}
\end{figure}
In~\cref{fig:3d_composition}, we illustrate how to compute all relevant sub-voxel volumes for vertically stacked detector elements.
For calculating the relative volumes of the individual voxel elements in a first step the cell is divided by a (magenta) plane that intersects the front and back plane into two halves.
Then the volumes are calculated from above and below.
For the upper and lower blue elements, the vectors \( \mathbf{C_{\mathbf{m}}} \) are \( (1, 1, 0, 1)^T \) and \( (0, 0, 1, 0)^T \), respectively.
To obtain the volume elements of the next (green) cells, they are subtracted from the subsequent cells with the vectors \( \mathbf{C_{\mathbf{m}}} = (1, 1, 1, 1)^T \), which reach the top and bottom surface.
Finally, the volumes of the divided middle cell must be added up.
% \begin{algorithm}[t]%
%     \caption{TODO Computation of the forward projection in three dimensions.}%
%     \label{algo:forward projection 3d}
%     \DontPrintSemicolon
%     \KwData{Image \( \mathbf{u} \) to be projected}
%     \For{pixels \( k = 1, 2, \dotsc, N\)}{
%         Compute relevant detectors \( \mathcal{D} \) with \cref{eq:relevant detectors 2d}.\;
%         \For{relevant detectors \( p \in \mathcal{D} \)}{
%             Compute angles \( \beta_1,\beta_2, \dotsc \) ... and depending on their relation compute the area factors .... and set \( \mathbf{p}_p = \) ...
%         }
%     }
%     \KwResult{Projections \( \mathbf{p} = \SystemMatrix \mathbf{u} \)}
% \end{algorithm}
% \begin{algorithm}[t]%
%     \caption{TODO Computation of the backprojection in three dimensions.}%
%     \label{algo:backprojection 3d}
%     \DontPrintSemicolon
%     \KwData{Measurements \( \Data \) to be backprojected}
%     \For{pixels \( k = 1, 2, \dotsc, N\)}{
%         Compute relevant detectors \( \mathcal{D} \) with \cref{eq:relevant detectors 2d}.\;
%         \For{relevant detectors \( p \in \mathcal{D} \)}{
%             Compute angles \( \beta \) ... and depending on their relation compute the area factors .... and set \( \mathbf{p}_p = \) ...
%         }
%     }
%     \KwResult{Backprojection \( \mathbf{u} = \SystemMatrix^\top \Data \)}
% \end{algorithm}
\section{Numerical Experiments}
In this section we show numerical results obtained from reconstructing images from sinogram data with the method as well as with the line-based reference implementation provided in the Astra toolbox~\cite{aarle2015astra,aarle2016astra}\footnote{Note that Astra also provides an area based forward operator for 2D, not, however, for 3D. We always use Astra's line-based operator for comparison in this paper.}.

For all experiments, given projection data $p$ we reconstruct the image of interest $u^\star$ via solving
\begin{equation}\label{eq:exp_tikh}
    \mathbf{u}^\star \in \argmin_{\mathbf{u} \in \R^N} F(u) \coloneqq \tfrac{1}{2}\|\SystemMatrix \mathbf{u} - \Data\|^2 + \tfrac{\lambda}{2}\| \mathbf{u} \|^2
\end{equation}
using \gls{nag}~\cite{nesterov269method}. If not specified otherwise, we iterate the algorithm until either $1000$ iterations are exceeded, or the squared Euclidean norm of the gradient of the $k$-th iterate $\mathbf{u}_k$ satisfies $\|\nabla F(\mathbf{u}_k)\|^2<\num{1e-9}$. The matrix $\SystemMatrix$ will always be either the line-based or the area-/volume-based \gls{sm}. In all cases we normalize the $\SystemMatrix$ by dividing by its spectral norm which is approximated using 100 power iterations. This way we can always choose the same parameter $\lambda$ for all methods\footnote{We noticed unstable behavior when using the Astra with the conjugate gradient method which we suspect to be a consequence of larger adjoint mismatch. Computing $|\langle \SystemMatrix\mathbf{u}, \mathbf{p}\rangle - \langle \mathbf{u},\SystemMatrix^* \mathbf{p}\rangle|$ with $\SystemMatrix,\SystemMatrix^*$ the provided forward and backward projections normalized by the spectral norm and random unit vectors $\mathbf{u},\mathbf{p}$ we obtain in 2D with our implementation $\num{1e-19}$ and Astra $\num{1e-10}$ and in 3D numerically zero with our method and $\num{1e-4}$ with Astra.}.

\subsection{The Two-Dimensional Case}
\subsubsection{Synthetic Data}
\begin{figure}
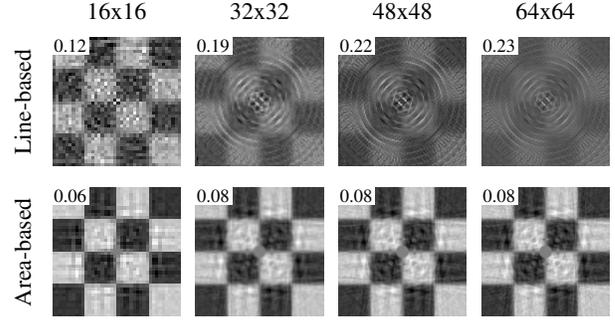

    \centering
    \begin{tikzpicture}[font=\small]
    \foreach [count=\ireso] \reso/\mseline/\msearea\resoim in {4x4/0.12/0.06/16x16, 8x8/0.19/0.08/32x32, 12x12/0.22/0.08/48x48, 16x16/0.23/0.08/64x64} {
        \node at (\ireso*1.9, -0.8) {\resoim};
            \begin{scope}
                \node (img) at (\ireso*1.9, -2)
                    {\includegraphics[width=1.7cm]{figures/2D-checkerboard/4x4/line-based/2D-\reso.png}};
                % Number label (adjust -0.65cm and +0.65cm if needed)
                \node[anchor=north west,
                    font=\scriptsize,
                    fill=white,
                    inner sep=1pt,
                    xshift = 0.1cm,
                    yshift=-0.1cm]
                    at (img.north west) {\mseline};
            \end{scope}
            \begin{scope}
                \node (img) at (\ireso*1.9, -4)
                    {\includegraphics[width=1.7cm]{figures/2D-checkerboard/4x4/area-based/2D-\reso.png}};
                % Number label (adjust -0.65cm and +0.65cm if needed)
                \node[anchor=north west,
                    font=\scriptsize,
                    fill=white,
                    inner sep=1pt,
                    xshift = 0.1cm,
                    yshift=-0.1cm]
                    at (img.north west) {\msearea};
            \end{scope}
    };
    \node [rotate=90] at (0.65, -2) {Line-based};
    \node [rotate=90] at (0.65, -4) {Area-based};
    \end{tikzpicture}
    \caption{%
        Reconstructed two-dimensional checkerboard images.
        Rows: Methods of system-matrix computation.
        Columns: Image resolution in pixels.
        In the top left corner of each image we show the MSE to the ground truth.
    }
    \label{fig:2D-checkerboard}
\end{figure}

We consider the task of reconstructing a 2D checkerboard pattern from its sinogram (\cf~\cref{fig:2D-checkerboard}) by solving the noiseless inverse problem. The projection data $\Data$ is generated as $\Data=\SystemMatrix \mathbf{u} + \mathbf{\epsilon}$ with $\SystemMatrix$ the area-based \gls{sm} and $\epsilon\sim \mathcal{N}(0,\sigma^2\mathrm{I})$ and $\sigma=10^{-4}$.
The parameters for this experiment are the number of projections $N_p = 60$, which are equidistant in $[0,2\pi)$, number of detector elements $N_d = 60$, $d_y=0.75$, $d?s=250$, and $a=b=1$. The regularization parameter is set to $\lambda=10^{-4}$.
The results are shown in~\cref{fig:2D-checkerboard}. In particular, there we successively increase the image resolution $N$, showing that the area-based method yields better results in the underdetermined regime than the line-based reconstruction.

\subsubsection{Real Data}
For an experiment using real \gls{ct} data we use the data set~\cite{hamalainen2015tomographic} which contains 2D \gls{ct} data of a walnut. The dataset includes downsampled sinograms with $N_p = 120$ and
$N_{d} \in \{82, 164, 328, 2296\}$ and the raw sinogram data with $N_p \times N_{d} = 1200 \times 2296$.
It also comes with a reference reconstruction that was generated as the \gls{fbp} from the high-resolution sinogram.
The results for $\lambda\in\{0.0001,0.01\}$, $N \in \{\num{200}\mathrm{x}\num{200},\num{500}\mathrm{x}\num{500}\}$ are shown in~\cref{fig:2D-walnut}. We find that the area based \gls{sm} provides the most benefit in regimes with little regularization.
\input{figures/2D-walnut.tex}
\subsubsection{Computational complexity}
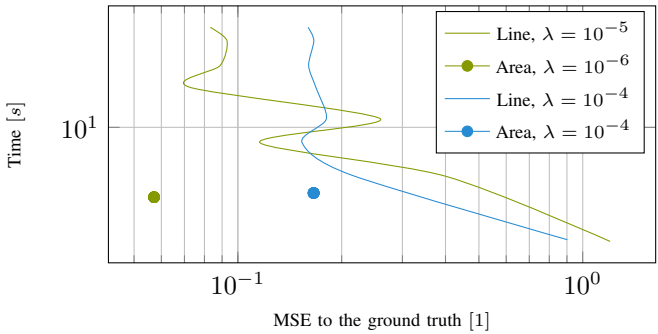
\begin{figure}
    \centering
\begin{tikzpicture}
  \begin{axis}[
    xlabel={\scriptsize MSE to the ground truth $[1]$},
    ylabel={\scriptsize Time $[s]$},
    grid=both,
    width=\linewidth,
    height=5cm,
    xmode=log,
    ymode=log,
  ]
  %   % Read x,y columns from a CSV file
  %   \addplot[
  %     color = plotred,
  %     % mark=*,
  %     smooth,
  %   ]
  %   table[
  %     col sep=comma,
  %     x=mse_line,
  %     y=time_line
  %   ] {figures/2D-complexity/lam_1e-06/complexity.txt};
  %   \addplot[
  %     color=plotred,
  %     mark=*,
  %     smooth,
  %   ]
  %   table[
  %     col sep=comma,
  %     x=mse_area,
  %     y=time_area
  %   ] {figures/2D-complexity/lam_1e-06/complexity.txt};

    \addplot[
      color = plotgreen,
      % mark=*,
      smooth,
    ]
    table[
      col sep=comma,
      x=mse_line,
      y=time_line
    ] {figures/2D-complexity/lam_1e-05/complexity.txt};
    \addplot[
      color=plotgreen,
      mark=*,
      smooth,
    ]
    table[
      col sep=comma,
      x=mse_area,
      y=time_area
    ] {figures/2D-complexity/lam_1e-05/complexity.txt};

    \addplot[
      color = plotblue,
      % mark=*,
      smooth,
    ]
    table[
      col sep=comma,
      x=mse_line,
      y=time_line
    ] {figures/2D-complexity/lam_0.0001/complexity.txt};
    \addplot[
      color=plotblue,
      mark=*,
      smooth,
    ]
    table[
      col sep=comma,
      x=mse_area,
      y=time_area
    ] {figures/2D-complexity/lam_0.0001/complexity.txt};
    \legend{{\scriptsize Line,
    % $\lambda=10^{-6}$}, {\scriptsize Area, $\lambda=10^{-6}$}, {\scriptsize Line, 
    $\lambda=10^{-5}$}, {\scriptsize Area, $\lambda=10^{-6}$}, {\scriptsize Line, $\lambda=10^{-4}$}, {\scriptsize Area, $\lambda=10^{-4}$},}
  \end{axis}
\end{tikzpicture}
\caption{Computational complexity of the proposed area-based and Astra's line-based reconstruction. For the line-based reconstruction we use multiple lines per detector element and compute their averages. As the number of lines per detector element increases, the reconstruction accuracy increases at the cost of longer computation times. We find no setting, where line-based reconstruction outperforms area based reconstruction in terms of accuracy \emph{and} speed.}
\label{fig:complexity}
\end{figure}

In~\cref{fig:complexity} we compare the computational complexity of reconstructing images with the area- and line-based \gls{sm}. More precisely, we generate a sinogram as again as $\Data = \SystemMatrix \mathbf{u}+\epsilon$ where now $\mathbf{u}$ is a 40x40 Shepp-Logan phantom and the noise level $\sigma=10^{-4}$. In \cref{fig:complexity} for different values of $\lambda$ we plot the computation times over the MSE to the ground truth image. In order to obtain multiple data points for the line-based method we compute the line-based projections with an increasing number of lines per detector element which are averaged afterwards. This way we approximate the area based discretization leading to decreasing reconstruction errors at the cost of more computational effort.
We find that for no setting the line based reconstruction is faster \emph{and} more accurate than the area based one. Moreover it seems that the increased accuracy of the line-based approximation saturates despite increasing the number of lines which is most likely due to the error caused by simply averaging multiple lines which is only correct in the limit in the parallel beam setting. Moreover, it should be noted that we could not use Astra's CUDA implementation for this experiment as it somehow led to decreasing accuracy when the number of lines per detector element was increased.
\subsection{The Three-dimensional Case}
In order to reduce computational complexity and simplify algorithms, in 3D we restrict to geometries where each image voxel contributes to at most three detector elements in the $z$ direction.
\subsubsection{Synthetic Data}
\begin{figure}
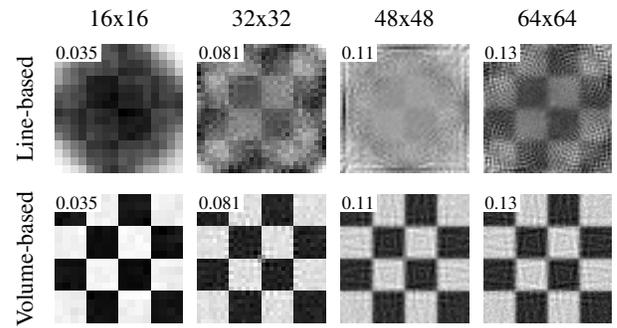

    \centering
    \begin{tikzpicture}[font=\small]
    \foreach [count=\ireso] \reso\resoname\mseline\msearea in {16x16/4x4/0.29/0.035, 32x32/8x8/0.41/0.081, 48x48/12x12/0.42/0.11, 64x64/16x16/0.42/0.13} {
        \node at (\ireso*1.9, -0.8) {\reso};
            \node (\ireso-line) at (\ireso*1.9, -1*2) {\includegraphics[width=1.7cm]{figures/3D-checkerboard/4x4/line-based/3D-\resoname_v3.png}};
            \node (\ireso-vol) at (\ireso*1.9, -2*2) {\includegraphics[width=1.7cm]{figures/3D-checkerboard/4x4/volume-based/3D-\resoname_v3.png}};
            \node[anchor=north west,
                    font=\scriptsize,
                    fill=white,
                    inner sep=1pt,
                    xshift = 0.1cm,
                    yshift=-0.1cm]
                    at (\ireso-line.north west) {\msearea};
            \node[anchor=north west,
                    font=\scriptsize,
                    fill=white,
                    inner sep=1pt,
                    xshift = 0.1cm,
                    yshift=-0.1cm]
                    at (\ireso-vol.north west) {\msearea};
        }
    \node [rotate=90] at (0.65, -2) {\small Line-based};
    \node [rotate=90] at (0.65, -4) {\small Volume-based};
    \end{tikzpicture}
    \caption{
        Second slice of the reconstructed three-dimensional checkerboard images.
        Rows: Methods of system-matrix computation.
        Columns: Image resolution in pixels.
    }
    \label{fig:3D-checkerboard-slices}
\end{figure}

Also in 3D we repeat the checkerboard experiment. In this case, the checkerboard consists of 4 black/white blocks in each direction ($x$, $y$, and $z$). We set the maximum number of iterations of the \gls{nag} algorithm to 500. All other parameters are kept as in the 2D experiment.
Again, we find that line-based discretization breaks down sooner in the underdetermined regime than the volume-based method, \cf~\cref{fig:3D-checkerboard-slices}.

\subsubsection{Real Data}
\begin{figure}
    \centering
    \begin{tikzpicture}[font=\small]
    \node at (1.9, -0.8) {Ground truth};
    \node (gt-line) at (1.9, -1*2) {\includegraphics[width=1.7cm]{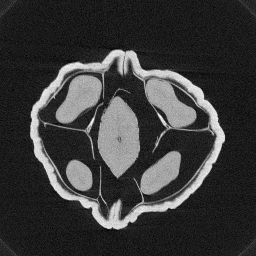}};
    \node (gt-vol) at (1.9, -2*2) {\includegraphics[width=1.7cm]{figures/3D-walnut/gt.png}};
    \foreach [count=\ilam] \lam\mseline\msearea in {0.001/0.059/0.048,0.005/0.038/0.040,0.1/0.040/0.036} {
        \node at (1.9+\ilam*1.9, -0.8) {$\lambda=\lam$};
            \node (\ilam-line) at (1.9+\ilam*1.9, -1*2) {\includegraphics[width=1.7cm]{figures/3D-walnut/line-based/lam_\lam.png}};
            \node (\ilam-vol) at (1.9+\ilam*1.9, -2*2) {\includegraphics[width=1.7cm]{figures/3D-walnut/volume-based/lam_\lam.png}};
            \node[anchor=north west,
                    font=\scriptsize,
                    fill=white,
                    inner sep=1pt,
                    xshift = 0.1cm,
                    yshift=-0.1cm]
                    at (\ilam-line.north west) {\mseline};
            \node[anchor=north west,
                    font=\scriptsize,
                    fill=white,
                    inner sep=1pt,
                    xshift = 0.1cm,
                    yshift=-0.1cm]
                    at (\ilam-vol.north west) {\msearea};
        }
    \node [rotate=90] at (0.65, -2) {\small Line-based};
    \node [rotate=90] at (0.65, -4) {\small Volume-based};
    \end{tikzpicture}
    \caption{
        Reconstruction of a 3D walnut. We show a slice at $z=10$. Columns: Ground truth image and results with different regularization parameters $\lambda$. Rows: Line- and volume-based reconstruction.
    }
    \label{fig:3D-walnut}
\end{figure}

We preform experiments using real \gls{ct} measurements from the data set \cite{meaney2022cone} which contains 3D \gls{ct} measurements of a walnut. 
The ground truth image is obtained via Astra's FDK implementation using the high resolution sinogram.
For the reconstructions we set $N_p=180$ evenly spaced in $[0,2\pi)$, $N_d = N_{d_y}\times N_{d_z} = 624\times 296$, $d_y=d_z=0.2$, $N=N_x\times N_y\times N_z=256\times 256\times 1$, $a=b=0.125$, $c=0.1$.
The results are shown in \cref{fig:3D-walnut}. In this experiment the visual differences between volume- and line-based reconstructions are negligible. However, on average we find slightly improved MSEs using the volume-based reconstruction. We believe that more research has to be conducted to determine in which parameter regimes the benefits are greatest in 3D.

\section{Conclusion}
We have presented a new method for the computation of a consistent \gls{sm} for cone beam flat detector \gls{ct}. In particular, the method does not rely on costly subroutines, but instead, we could derive explicit formulas for the entries of the \gls{sm}. We have shown in experimental results in 2D and 3D that the consistent \gls{sm} leads to improved reconstruction quality compared to conventional line-based methods which is particularly relevant for quantitative \gls{ct}. Moreover, the \gls{sm} is more stable in the underdetermined regime potentially allowing for improvements in low-dose and dynamic \gls{ct}.
\subsection{Limitations and future work}
While improvements in terms of reconstruction quality are evident, a limitation of the proposed method is still the computational burden for higher resolutions. In future work we plan to improve upon the current implementation to further facilitate practical adoption of the method.
Additionally, we plan to investigate the empirical behavior of the area-/volume-based discretization in comparison to line-based methods in more detail. We also aim to extend the proposed strategy of explicitly computing volume based weighting coefficients to other geometries beyond flat detector \gls{ct} and to measurement during rotation and longitudinal translation as in spiral \gls{ct}. 
\printbibliography

\appendix

\subsection{Postponed derivations}
\subsubsection{Fundamental relations}\label{sec:appendix}
% \subsubsection{Fundamental relations}\label{sec:fundamentals}
We begin by deriving various fundamental relations.
Defining $\beta_{\mathbf{m}_y}\in (-\pi/2,\pi/2)$ via 
\begin{equation}\label{eq:beta_computation}
    \tan\beta_{\mathbf{m}_y} = \frac{\mathbf{m}_{y}d_{y}}{s + d}
\end{equation}
We can equivalently write $\dvec_{\mathbf{m}}$ from \eqref{eq:z rotation} as
\begin{equation}
	\dvec_{\mathbf{m}}  = \begin{pmatrix}
		s\cos\varphi \\
		s\sin\varphi \\
		0
	\end{pmatrix} + \lambda \begin{pmatrix}
		L_\mathbf{m}\cos(\varphi - \beta_{\mathbf{m}_y}) \\
		L_\mathbf{m}\sin(\varphi - \beta_{\mathbf{m}_y}) \\
		 - \mathbf{m}_{z}d_{z}
	 \end{pmatrix}
\end{equation}
with
\begin{equation}
    L_{\mathbf{m}} = \sqrt{{(s + d)}^{2} + \mathbf{m}_{y}^{2}d_{y}^{2}} = (s+d)/\cos\beta_{\mathbf{m}_y}. 
\end{equation}
% Thus, the ray \( \dvec_{\svec} \) and a vector from $\svec$ to the origin enclose the angle \( \beta_{m} \).
It can be concluded that for an arbitrary point \( (x, y, z) \) along a trajectory with the angle \( \beta \), the following relations are valid:
\begin{equation}
\begin{split}
	z &= \frac{\mathbf{m}_{z}d_{z}}{s+d}\frac{\cos\beta}{\cos{(\varphi} - \beta)}\left( s\cos\varphi - x \right)\\
    &= \frac{\mathbf{m}_{z}d_{z}}{s+d}\frac{\cos\beta}{\sin {(\varphi} - \beta)}\left( s\sin\varphi - y \right),
 \label{eq:z expr 1}
 \end{split}
\end{equation}
\begin{equation}
\begin{split}
	z &= \left\lbrack \left( s\cos\varphi - x \right)\cos\varphi + \left( s\sin\varphi - y \right)\sin\varphi \right\rbrack\frac{\mathbf{m}_{z}d_{z}}{s + d}\\
    &= \left( s - x\cos\varphi - y\sin\varphi \right)\frac{\mathbf{m}_{z}d_{z}}{s + d},
\end{split} 
\end{equation}
\begin{equation}
\begin{split}
	x &= s\cos\varphi - \frac{\cos(\varphi- \beta)}{\cos\beta}\frac{s+d}{\mathbf{m}_{z}d_{z}}z \\
    &= s\cos\varphi + \frac{1}{\cos\varphi}\left\lbrack \left( s\sin\varphi - y \right)\sin\varphi - \frac{s+d}{\mathbf{m}_{z}d_{z}}z \right\rbrack,
 \label{eq:x expr 1}
\end{split} 
\end{equation}
\begin{equation}
\begin{split}
	y &= s\sin\varphi - \frac{\sin(\varphi- \beta)}{\cos\beta}\frac{s+d}{\mathbf{m}_{z}d_{z}}z \\
    & = s\sin\varphi + \frac{1}{\sin\varphi}\left\lbrack \left( s\cos\varphi - x \right)\cos\varphi - \frac{s+d}{\mathbf{m}_{z}d_{z}}z \right\rbrack,
  \label{eq:y expr 1}
\end{split}  
\end{equation}
\begin{equation}
\begin{split}
    x &= - s\frac{\sin\beta}{\sin{(\varphi} - \beta)} + y\cot{(\varphi} - \beta)\\ 
    &= s\cos\varphi - (s\sin\varphi - y)\cot{(\varphi} - \beta),
    \label{eq:x expression}
\end{split}   
\end{equation}

\begin{equation}
\begin{split}
    y &= s\frac{\sin\beta}{\cos{(\varphi} - \beta)} + x\tan
    {(\varphi} - \beta) \\ 
    &= s\sin\varphi - (s\cos\varphi - x)\tan{(\varphi} - \beta),
    \label{eq:y expression}
\end{split}    
\end{equation}

\begin{equation}
    \tan{(\varphi} - \beta) = \frac{s\sin\varphi - y}{s\cos\varphi - x}.
    \label{eq:tan phi minus beta}
\end{equation}
From~\cref{eq:x expression}, we may derive that
\begin{equation}
    \tan\beta = \frac{y\cos\varphi - x\sin\varphi}{s - y\sin\varphi - x\cos\varphi }.
    \label{eq:tan beta}
\end{equation}
% \subsection{The Two-Dimensional Case}
% Via these $\pi/2$-rotations the $x$ and $y$ coordinates $x_1<x_2,y_1<y_2$ of a pixel $[x_1,x_2]\times[y_1,y_2]$ swap according to
% For this purpose, the fourfold symmetry of rectangular or quadratic pixels with respect to the 360° circumference is used:
% \begin{equation}
%     \begin{aligned}
%         &[x_1,x_2,y_1,y_2]\rightarrow
%         [y_1,y_2,-x_2,-x_1]\rightarrow\\
%         &[-x_2,-x_1,-y_2,-y_1]\rightarrow
%         [-y_2,-y_1,x_1,x_2]
		% &[n_xa, (n_x + 1)a, n_yb, (n_y+1)b] \rightarrow \\& [n_yb, (n_y+1)b, -n_xa, -(n_x+1)a] \rightarrow \\
		% &[-(n_x+1)a, -n_xa, -(n_y+1)b, -n_yb] \rightarrow \\& [-(n_y+1)b, -n_yb, n_xa, (n_x+1)a]
% 	\end{aligned}
% \end{equation}
% Note that in the second and fourth case the sides \( a \) and \( b \) be must be interchanged.
From~\cref{eq:tan beta}, the following relations can be derived:
\begin{equation}\label{eq:gamma_computation}
\begin{split}
	\tan\gamma_{\mathbf{n},2,\phi} &= \frac{\mathbf{n}_{y}b\cos\varphi - \left( \mathbf{n}_{x} + 1 \right)a\sin\varphi}{s- \mathbf{n}_{y}b\sin\varphi - \left( \mathbf{n}_{x} + 1 \right)a\cos\varphi},\\
	\tan\gamma_{\mathbf{n},1,\phi} &= \frac{\mathbf{n}_{y}b\cos\varphi - \mathbf{n}_{x}a\sin\varphi}{s - \mathbf{n}_{y}b\sin\varphi - \mathbf{n}_{x}a\cos\varphi},\\
	\tan\gamma_{\mathbf{n},4,\phi} &= \frac{{(\mathbf{n}}_{y} + 1)b\cos\varphi - \mathbf{n}_{x}a\sin\varphi}{s - {(\mathbf{n}}_{y} + 1)b\sin\varphi - \mathbf{n}_{x}a\cos\varphi},\text{ and}\\
	\tan\gamma_{\mathbf{n},3,\phi} &= \frac{{(\mathbf{n}}_{y} + 1)b\cos\varphi - \left( \mathbf{n}_{x} + 1 \right)a\sin\varphi}{s - {(\mathbf{n}}_{y} + 1)b\sin\varphi - \left( \mathbf{n}_{x} + 1 \right)a\cos\varphi}.
    \end{split}
\end{equation}
From~\cref{eq:tan phi minus beta}, it follows
\begin{equation}
\begin{split}
    \cot (\varphi - \gamma_{\mathbf{n},1,\phi}) - \cot (\varphi - \gamma_{\mathbf{n},2,\phi})) 
    % \\= \frac{s\cos\varphi - n_{x}a}{s\sin\varphi - n_{y}b} - \frac{s\cos\varphi - a\left( n_{x} + 1 \right)}{s\sin\varphi - n_{y}b} \\
    = \frac{a}{s\sin\varphi - \mathbf{n}_{y}b}.
 \label{eq:cot diff}
\end{split} 
\end{equation}
Thus, the following relation results,
\begin{equation}
	s\sin\varphi - \mathbf{n}_{y}b = a\frac{\sin{(\varphi - \gamma_{\mathbf{n},2,\phi})}\sin{(\varphi - \gamma_{\mathbf{n},1,\phi})}}{\sin(\gamma_{\mathbf{n},1,\phi} - \gamma_{\mathbf{n},2,\phi})}.
 \label{eq:sin phi ny b}
\end{equation}
Combining~\cref{eq:tan phi minus beta} and~\cref{eq:sin phi ny b}, we obtain
\begin{equation}
	s\cos\varphi - \mathbf{n}_{x}a = a\frac{\sin{(\varphi - \gamma_{\mathbf{n},2,\phi})}\cos{(\varphi - \gamma_{\mathbf{n},1,\phi})}}{\sin(\gamma_{\mathbf{n},1,\phi} - \gamma_{\mathbf{n},2,\phi})},
  \label{eq:cos phi nx a}
\end{equation}
\begin{equation}
	s\cos\varphi - \left( \mathbf{n}_{x} + 1 \right)a = a\frac{\cos{(\varphi - \gamma_{\mathbf{n},s,\phi})}\sin{(\varphi - \gamma_{\mathbf{n},1,\phi})}}{\sin(\gamma_{\mathbf{n},1,\phi} - \gamma_{\mathbf{n},2,\phi})}.
 \label{eq:cos phi nx1 a}
\end{equation}
% The angle \( \beta_{\mathbf{m}} \) of the ray from \( \svec \) to the \( \mathbf{m} \)-th intersection of the detector elements can be calculated from
% \begin{equation}
% 	\tan\beta_{\mathbf{m}} = \mathbf{m}_{y}d_{y}/(s+d).
% \end{equation}
% \subsubsection{Triangle Area}
% Using~\cref{eq:sin phi ny b} and~\cref{eq:cos phi nx1 a}, the two sides of the obtained rectangular triangle can be calculated as follows:
% \begin{equation}
% \begin{split}
% 	&\left| \mathbf{p}_{\mathbf{x}} - \mathbf{v}_{1} \right|  \\
%     =&\left|s\frac{\sin{\beta_{i}}}{\cos{(\varphi - \beta_{i})}}-n_yb+(n_x+1)a\tan{(\varphi - \beta_{i})}\right|
%     \\ = &a\left| \frac{\sin{(\varphi} - \betagtwo)\sin{(\beta_{i}} - \betago)}{\sin{(\betagtwo} - \betago)\cos{(\varphi - \beta_{i})}} \right|
% \end{split}
%  \label{eq:px v1}
% \end{equation}
% and
% \begin{equation}
% \begin{split}
% 	&\left| \mathbf{v}_{1}{- \mathbf{p}}_{\mathbf{y}} \right| \\
%     =&\left|s\frac{\sin{\beta_{i}}}{\sin{(\varphi - \beta_{i})}}-n_yb\cot{(\varphi - \beta_{i})}+(n_x+1)a\right|
%     \\ = &a\left| \frac{\sin{(\varphi} - \betagtwo)\sin{(\beta_{i}} - \betago)}{\sin{(\betagtwo} - \betago)\sin{(\varphi - \beta_{i})}} \right|.
%   \label{eq:v1 py}
% \end{split}
% \end{equation}
% Hence, the resulting area \( A_{\text{tri}} \) can be written as:
% \begin{equation}
% \begin{split}
% 	&A_{\text{tri}}(\varphi,\betago,\betagtwo,\beta_{i})\\
%     = &{\frac{a^{2}}{\left| \sin{(2\varphi - {2\beta}_{i})} \right|}\left( \frac{\sin{(\varphi} - \betagtwo)\sin{(\beta_{i}} - \betago)}{\sin{(\betagtwo} - \betago)} \right)}^{2}.
% \end{split}
% \end{equation}

\subsubsection{The Three-dimensional Case}
\paragraph{Prismatoid}
The volume of the prismatoid(\cf~\cref{fig:3d_cases} (c)) is computed as
\begin{equation}
	\begin{aligned}
		\vol = \frac{a}{6}\left\lbrack 2\frac{\left| V_{1}P_{1} \right| + \left| W_{1}Q_{1} \right|}{2}\left| V_{1}W_{1} \right| + \right. \\
        + \frac{\left| V_{1}P_{1} \right| + \left| W_{1}Q_{1} \right|}{2}\left| V_{0}W_{0} \right| +\\
		+ \frac{\left| V_{0}P_{0} \right| + \left| W_{0}Q_{0} \right|}{2}\left| V_{1}W_{1} \right| + \\
        \left. + 2\frac{\left| V_{0}P_{0} \right| + \left| W_{0}Q_{0} \right|}{2}\left| V_{0}W_{0} \right| \right\rbrack
	\end{aligned}
\end{equation}
leading to
\begin{equation}
\begin{split}
	\vol = \frac{1}{6}a^{3}\ \left| \frac{\tan\alpha_{\mathbf{m}_{z}}}{\sin\varphi}\lbrack \left( f_{\mathbf{m}_y}^{2} + f_{\mathbf{m}_y}g_{\mathbf{m}_y} + g_{\mathbf{m}_y}^{2} \right) - \right.\\
    \left.- \left( f_{\mathbf{m}_y + 1}^{2} + f_{\mathbf{m}_y + 1}g_{\mathbf{m}_y + 1} + g_{\mathbf{m}_y + 1}^{2} \right) \rbrack \ \right|.
 \end{split}
\end{equation}
Using the relation~\cref{eq:f minus g}, the final result is:
\begin{equation}
\begin{split}
	\vol = &\frac{1}{6}a^{3}\bigg\lvert \frac{\tan\alpha_{\mathbf{m}_{z}}}{\sin\varphi}\left\lbrack \frac{\cos{(\varphi-\beta_{\mathbf{m}_y})}}{\cos{\beta_{\mathbf{m}_y}}}\left( g_{\mathbf{m}_y}^{3} - f_{\mathbf{m}_y}^{3} \right) - \right. \\
    &\left. - \frac{\cos{(\varphi-\beta_{\mathbf{m}_y+1})}}{\cos{\beta_{\mathbf{m}_y+1}}}{\left( g_{\mathbf{m}_y + 1}^{3} - f_{\mathbf{m}_y + 1}^{3} \right)}\right\rbrack \bigg\rvert
\end{split}
\end{equation}

\subsubsection{Top-view triangle}\label{sec:appendix chopped prism}
In the case~\cref{fig:3d_cases} (d) the volume is computed as
\begin{equation}
	\vol = \frac{1}{6}\left| V_{12}W_{1} \right|\left| V_{12}W_{2} \right|\left\lbrack \left| Q_{1}W_{1} \right| + |PV_{12}| + \left| Q_{2}W_{2} \right| \right\rbrack
\end{equation}
where
\begin{equation}
\begin{split}
	\left| V_{12}W_{1} \right| &= \frac{a}{\sin\varphi}\left\lbrack g_{\mathbf{m}_y} - g_{g}\right\rbrack\\ 
    &= a\frac{\sin{(\varphi} - \gamma_{\mathbf{n},1,\phi})\sin{(\beta_{\mathbf{m}_y}} - \gamma_{\mathbf{n},2,\phi})}{\sin{(\gamma_{\mathbf{n},1,\phi}} - \gamma_{\mathbf{n},2,\phi})\cos{(\varphi - \beta_{\mathbf{m}_y})}},
 \label{eq:v11 w1}
 \end{split}
\end{equation}
\begin{equation} 
\begin{split}
	\left| V_{12}W_{2} \right| &= \frac{a}{\cos\varphi}\left\lbrack g_{g} - h_{\mathbf{m}_y} \right\rbrack \\
    &= a\frac{\sin{(\varphi} - \gamma_{\mathbf{n},1,\phi})\sin{(\beta_{\mathbf{m}_y}} - \gamma_{\mathbf{n},2,\phi})}{\sin{(\gamma_{\mathbf{n},1,\phi}} - \gamma_{\mathbf{n},2,\phi})\sin{(\varphi - \beta_{\mathbf{m}_y})}},
 \label{eq:v11 w2}
 \end{split}
\end{equation}
\begin{equation}
	\left| Q_{1}W_{1} \right| = ag_{\mathbf{m}_y}\tan\alpha_{\mathbf{m}_{z}},
\end{equation}
\begin{equation}
	\left| PV_{12} \right| = ag_{g}\tan\alpha_{\mathbf{m}_{z}}
 \label{eq:p v11}
\end{equation}
and
\begin{equation}
	\left| Q_{2}W_{2} \right| = ah_{\mathbf{m}_y}\tan\alpha_{\mathbf{m}_{z}},
\end{equation}
with
\begin{equation}
\begin{split}
	g_{g} &= \frac{c}{a}\frac{\mathbf{n}_{z} + 1}{\tan\alpha_{\mathbf{m}_{z}}} - \frac{s - \left( \mathbf{n}_{x} + 1 \right)a\cos\varphi - \mathbf{n}_{y}b\sin\varphi}{a} \\ 
    &= \frac{c}{a}\frac{\mathbf{n}_{z} + 1}{\tan\alpha_{\mathbf{m}_{z}}} - \frac{\cos\gamma_{\mathbf{n},2,\phi}\sin{(\varphi} - \gamma_{\mathbf{n},1,\phi})}{\sin{(\gamma_{\mathbf{n},1,\phi}} - \gamma_{\mathbf{n},2,\phi})},
\end{split}
\end{equation}
\begin{equation}
\begin{split}
	h_{\mathbf{m}_y} &= \frac{c}{a}\frac{\mathbf{n}_{z} + 1}{\tan\alpha_{\mathbf{m}_{z}}}\ - \frac{\cos\beta_{\mathbf{m}_y}}{\sin{(\varphi} - \beta_{\mathbf{m}_y})}\frac{s\sin\varphi - \mathbf{n}_{y}a}{a}\\ 
    &= \frac{c}{a}\frac{\mathbf{n}_{z} + 1}{\tan\alpha_{\mathbf{m}_{z}}} - \frac{\cos\beta_{\mathbf{m}_y}\sin{(\varphi} - \gamma_{\mathbf{n},2,\phi})\sin{(\varphi} - \gamma_{\mathbf{n},1,\phi})}{\sin(\varphi - \beta_{\mathbf{m}_y})\sin(\gamma_{\mathbf{n},1,\phi} - \gamma_{\mathbf{n},2,\phi})}.
\end{split}
\end{equation}
The parameters \( g_{g} \), \( g_{\mathbf{m}_y} \) and \( h_{\mathbf{m}_y} \) obey the relation
\begin{equation}
	\cos\varphi\frac{\cos{(\varphi-\beta_{\mathbf{m}_y})}}{\cos{\beta_{\mathbf{m}_y}}} g_{\mathbf{m}_y} + \sin\varphi \frac{\sin{(\varphi-\beta_{\mathbf{m}_y})}}{\cos{\beta_{\mathbf{m}_y}}}h_{\mathbf{m}_y} = g_{g}.
 \label{eq:gi gg}
\end{equation}
% Note that~\cref{eq:v11 w1} and~\cref{eq:v11 w2} correspond to~\cref{eq:px v1} and~\cref{eq:v1 py}, which are the sides of the two-dimensional triangle.
After some algebraic work and using~\cref{eq:gi gg}, it is obtained:
\begin{equation}
\begin{split}
	\vol = &a^{3}\left| \frac{\tan\alpha_{\mathbf{m}_{z}}}{\sin\varphi}\left\lbrack\frac{\cos{(\varphi-\beta_{i})}}{\cos{\beta_{\mathbf{m}_y}}} g_{\mathbf{m}_y}^{3} + \right. \right.
    \\ &\left.\left. + \tan\varphi\frac{\sin{(\varphi-\beta_{i})}}{\cos{\beta_{\mathbf{m}_y}}} h_{\mathbf{m}_y}^{3}- \frac{1}{cos\varphi}g_{g}^{3}\right\rbrack \right|.
\end{split}
\end{equation}
\end{document}